\documentclass{article}

    \PassOptionsToPackage{numbers, compress}{natbib}
\usepackage[square,sort,comma,numbers]{natbib}

\usepackage{amsthm}
\usepackage{graphicx}
\usepackage[final]{neurips_2022}

\usepackage{amsmath,amsfonts,amssymb}
\usepackage{algorithm}
\usepackage{algorithmic}
\usepackage{dsfont}
\usepackage{xcolor} 
\usepackage{caption,color}
\usepackage{subcaption}
\usepackage{hyperref}
\usepackage{ulem}% colors
\usepackage[T1]{fontenc}
\usepackage[utf8]{inputenc}
\usepackage{enumitem}
\usepackage{comment}
\newcommand{\ust}{^{\star}}

\newcommand{\te}{\theta}

\newcommand{\cE}{\mathcal{E}}

\newcommand{\alig}[1]{\begin{align}#1\end{align}}
\newtheorem{assumption}{\textbf{Assumption}}

\newtheorem{theorem}{Theorem}[section]

\newtheorem{lemma}[theorem]{Lemma}
\newtheorem{remark}{Remark}
\title{Augmented RBMLE-UCB Approach for Adaptive Control of Linear Quadratic Systems}

\author{%
  Akshay Mete \\
  Texas A \& M University\\
  College Station, Texas, USA\\
  \texttt{akshaymete@tamu.edu} \\
  \And
   Rahul Singh \\
Indian Institute of Science \\
   Bengaluru, Karnataka, India \\
   \texttt{rahulsingh@iisc.ac.in} \\
   \AND
  P. R. Kumar\\
  Texas A \& M University\\
  College Station, Texas, USA\\
  \texttt{prk@tamu.edu} \\
}

\begin{document}

\maketitle
\begin{abstract}
We consider the problem of controlling an unknown stochastic linear system with  quadratic costs -- called the adaptive LQ control problem. We re-examine an approach called ``Reward-Biased Maximum Likelihood Estimate'' (RBMLE) \cite{kumar1982new} that was proposed more than forty years ago, and which predates the ``Upper Confidence Bound'' (UCB) method as well as the definition of ``regret'' for bandit problems \cite{lai1985asymptotically}. It simply added a term favoring parameters with larger rewards to the criterion for parameter estimation.  We show how the RBMLE and UCB methods can be reconciled, and thereby propose an Augmented RBMLE-UCB algorithm that combines the penalty of the RBMLE method with the constraints of the UCB method \cite{abbasi_11}, uniting the two approaches to optimism in the face of uncertainty. We establish that theoretically, this method retains $\Tilde{\mathcal{O}}(\sqrt{T})$ regret, the best known so far. We further compare the empirical performance of the proposed Augmented RBMLE-UCB and the standard RBMLE (without the augmentation) with UCB, Thompson Sampling, Input Perturbation, Randomized Certainty Equivalence and StabL on many real-world examples including flight control of Boeing 747 and Unmanned Aerial Vehicle. We perform extensive simulation studies showing that the Augmented RBMLE consistently outperforms UCB, Thompson Sampling and StabL by a huge margin, while it is marginally better than Input Perturbation and moderately better than Randomized Certainty Equivalence. 
% Unlike OFU and TS, Augmented RBMLE is able to ensure system stability even during the ``transient phase.'' 
\end{abstract}
\section{Introduction}
\label{sec:intro}
Consider a linear stochastic system $x_{t+1}=A\ust x_t+B\ust u_t+w_{t+1}$, where $x_{t}\in\mathds{R}^{n}$, $u_{t}\in \mathds{R}^m$ and
$w_t \in \mathds{R}^{n}$ are the state, control applied, and state noise, respectively, at time $t$.
We study the adaptive control\slash reinforcement learning problem \cite{aastrom2013adaptive,kumar1985survey,kumar1986stochastic,sutton2018reinforcement} where the controller's goal is to minimize the expected finite horizon quadratic cost 
$\mathds{E} \sum_{t=1}^{T} \left( x_t^{\intercal}Qx_t + u^{\intercal}_t Ru_t \right)$ by choosing the control input based on observing the past states and controls, without the knowledge of  ``true parameter'' $\te\ust := (A\ust,B\ust)^{\intercal}$. The deviation from what would have been optimally possible had the true parameter been known is measured by the ``regret,'' \cite{lai1985asymptotically},
% The goal of the agent is to choose the controls $\{u_t\}$ , so as to minimize the expected value of the cumulative cost incurred during $T$ steps, i.e. 
% This objective serves the purpose of keeping the state of the linear system close to the origin by using a minimal amount of control cost, and is a very popular approach to model real-world engineering control systems~
%  When the system matrices $(A,B)$ are known, and $w_t$ are independent, identically distributed across time and Gaussian, this is called the Linear Quadratic Gaussian (LQG) control problem.~In this work we study the \textit{adaptive LQG} problem, in which the controller\slash agent aims to minimize the LQG cost without knowing the parameter $\theta\ust := (A,B)$. More precisely, the agent's goal is to choose $u_t$ adaptively on the basis of the information collected during the system operation, so as to minimize the expected value of the ``regret''~\cite{lattimore2020bandit}
% \begin{equation*} \label{A,B linear system}
     $R(T) := \sum_{t=1}^{T} \left( x^{\intercal}_tQx_t + u^{\intercal}_t Ru_t \right) -T~J\ust(\theta\ust)$,
% \end{equation*}
where $J\ust(\te\ust):=
\limsup_{T \to \infty} \frac{1}{T} \mathds{E} \sum_{t=1}^{T} \left( x^{\intercal}_t Qx_t + u^{\intercal}_t Ru_t \right)$ is the optimal average cost achievable when $(A\ust,B\ust)$ are known.

%\prk{Do you use $T J(\theta)$ or the optimal cost over $[0,T]$?
%The former does not depend on the initial condition, whereas the latter does.}
%\am{We are using above regret equation for simulations plots.()}
% This is the adaptive LQG control problem that we address.
% Let $\pi^{\te}(\cdot)$ be an optimal stationary control policy~\cite{puterman2014markov} when the parameter $\te$ is known. 
For this adaptive linear-quadratic (LQ) control problem, 
% or for the problem of Adaptive Control of Markov Decision Processes (MDPs), 
there are broadly four classical approaches: the Reward-Biased Maximum Likelihood Estimate (RBMLE) approach \cite{kumar1982new}, the Diminishing Excitation (DE) approach \cite{chen1991identification}, the Upper Confidence Bound (UCB) approach \cite{lai1985asymptotically}, 
and the Thompson Sampling (TS) \cite{thompson1933likelihood} approach based on sampling from the posterior distribution. The RBMLE and UCB approaches are ``certainty equivalent'' (CE) in the sense that they  make an estimate ${\hat{\te}}_t$ of the unknown true parameter, and then take an action $u_t=K({\hat{\te}}_t)x_t$ that would be optimal if the estimate ${\hat{\te}}_t$ were indeed the true parameter. They only differ in what parameter estimate  they choose. The DE approach applies  $u_t=K({\hat{\te}}_t) x_t+v_t$, where $v_t$ is an added ``excitation,'' an independent noise, of diminishing variance and $\hat{\te}_t$ is the maximum likelihood estimate. The Randomized Certainty Equivalence (RCE)~\cite{faradonbeh_rce} adds excitation to the parameter estimate.

% In this paper we reconcile the RBMLE and UCB approaches, and armed with that, show how one may coherently combine the two approaches to design
% an Augmented RBMLE (ARBMLE) scheme. This is reminiscent of the Augmented Lagrangian approach
% to constrained optimization \cite{hestenes1969multiplier,powell1969method} that combines the primal penalty approach with a dual approach
% to constrained optimization. 

\subsection{The Contributions}
\begin{enumerate}[label=\arabic*.]
    \item We unite the two approaches to ``optimism under uncertainty'',
    RBMLE \cite{kumar1982new} and UCB \cite{lai1985asymptotically}, by showing that the RBMLE method is a penalty version of the constrained optimization problem of UCB for the case of linear quadratic systems \cite{kumar1983optimal,campi1998adaptive,prandini2000adaptive,bittanti2006adaptive,abbasi_11}. Based on this we propose an Augmented RBMLE (ARBMLE) method that combines the penalty and constrained versions 
    so that on the one hand it retains the analytical tractability of UCB and on the other hand provides the performance of RBMLE.
% Certainty equivalent schemes suffer from the closed-loop identifiability problem 
% \cite{mandl1974estimation,borkar1979adaptive} wherein the
% limiting parameter estimate identify the closed-loop behavior but not the open loop behavior.
% It was noticed in \cite{kumar1982new} that this implies that the limiting parameter estimate
% is biased in that the optimal cost associated with the model it parametrizes is necessarily larger than the true parameter. Motivated by this it was proposed in \cite{kumar1982new} to bias the parameters in the reverse way in favor of
% parameters associated with lower optimal costs -- the RBMLE method. This represented an optimistic approach favoring parameters with lower costs. Subsequently a different approach was proposed in \cite{lai1985asymptotically} which directly favored parameers with lower costs that could plausibly explain the data -- the UCB method. Here we unite these approaches to optimism in the face of uncertainty by showing that the RBMLE method is  a penalty version of an constrained optimization problem represented by UCB. Based on this observation, we unite the two methods by proposing an Augmented RBMLE (ARBMLE) that combines 
% This therefore optimal cost for the limiting estimate
% The RBMLE method \cite{kumar1982new} was motivated by the problem of overcoming the    where the parameter estimators can asymptotically identify the closed-loop behavior converge to a limit $\theta_\infty$ scan asymptotically identify the cl
    \item 
    We determine how to choose the biasing factor $\alpha(t)$ for ARBMLE, and establish
% resolve several critical questions to tailor RBMLE to LQG control.
% We determine 
a finite time regret bound 
$\Tilde{\mathcal{O}}(\sqrt{T})$,
% for the Augmented RBMLE-UCB algorithm for adaptive LQG,
the same as the OFULQ algorithm of~\cite{abbasi_11}, the best order available to date. 
% We analyze the regret of the ARBMLE method for adaptive LQ control, and establish that with a high probability it can be bounded as $\tilde{O}(\sqrt{T})$, the best known so far. 
   \item We perform extensive comparative simulation studies of the performance of:
   \begin{enumerate}
        \item ARBMLE and the standard RBMLE.
       \item OFULQ \cite{abbasi_11}, which is the UCB-approach adapted to the LQ problem.
       \item TS \cite{abeille_17} which is the Thompson sampling approach adapted to the LQ problem
       \item Input Perturbation (IP) \cite{faradonbeh2018input} which is a recent reincarnation of DE for the adaptive LQ problem that additionally assumes apriori knowledge of a stabilizing controller.
       \item Stabl \cite{lale22}, a modified OFULQ that adds ``excitation'' to the input for initialization.
       \item RCE \cite{faradonbeh_rce}, which adds excitation to the Least Squares Estimate (LSE).
   \end{enumerate}
   The examples used for our simulation study have been used in many recent papers \cite{dean_17,abbasi_19,lale22}, namely (a) the longitudinal flight control of Boeing 747 with linearized dynamics \cite{lale22},(b) Unmanned Aerial Vehicle (UAV) \cite{Zhao21,lale22} (c) unstable Laplacian dynamics \cite{dean_17}, and (d) large transient dynamics \cite{dean_17}. Our simulation results show that ARBMLE outperforms OFULQ, TS and StabL by a large margin, which is primarily due to lack of stabilization experienced by OFULQ and TS \cite{lale22}. RCE also exhibits a higher regret than ARBMLE. While the empirical performance of IP is marginally worse than ARBMLE. The results show that the ARBMLE has the same performance as the original RBMLE \cite{campi_98}, and they both outperform all the above other algorithms. Notably the choices made by ARBMLE and OFULQ within the confidence interval are very far apart, with ARBMLE outperforming OFULQ by a large margin. 
  
%   In this paper, we determine how to choose the biasing factor $\alpha(t)$, and establish
% % resolve several critical questions to tailor RBMLE to LQG control.
% % We determine 
% a finite time regret bound 
% $\Tilde{\mathcal{O}}(\sqrt{T})$
% for the Augmented RBMLE-UCB algorithm for adaptive LQG,
% which is the same as the OFU algorithm of~\cite{abbasi_11}, the best order available to date. 
%     We also comparatively evaluate the empirical performance of the Augmented RBMLE algorithm in extensive experiments.\\

% The Augmented RBMLE-UCB scheme not only has a $\Tilde{\mathcal{O}}(\sqrt{T})$ order of regret, which matches the best order of regret known to date \cite{abbasi_11}, \textcolor{violet}{but also in our simulation studies had considerably better regret than either the OFU algorithm \cite{abbasi_11}, which follows the UCB approach, or the TS \cite{abeille_17} algorithm.}

\end{enumerate}

\subsection{Previous Works}
Prior work on RBMLE has concentrated on establishing long-term average optimality \cite{kumar1982adaptive,kumar1982optimal,kumar1982new,kumar1983optimal,kumar1983simultaneous,borkar1990kumar,campi1998adaptive,prandini2000adaptive,bittanti2006adaptive,stettner1993nearly,borkar1991self,duncan1994almost}. %\rs{ RBMLE %was also shown to be stable in adaptive LQG %control~\cite{campi_98,prandini2000adaptive,bittanti2006adaptive,campi1998adaptive}, i.e. %$\limsup_{T\to\infty}\frac{1}{T}\sum_{t=1}^{T}\|x(t)\|^{2}<\infty$ almost surely.
%}
% \prk{Cite all of them in one shot}
    Recently, its regret performance has been addressed for Multi-Armed Bandits (MABs) \cite{liu2020exploration},
linear contextual bandits \cite{hung2020reward}, and
MDPs \cite{mete2021reward}. %\prk{ADD THE WORK OF \cite{campi1998adaptive,prandini2000adaptive} TPP}
% The seminal work on adaptive LQG control~\cite{kumar1983optimal} showed that the RBMLE algorithm attains the optimal infinite horizon average cost asymptotically by utilizing reward-biased estimates of the system parameters in order to generate controls.~However, it assumed that the unknown parameters can assume values from only a discrete set. This assumption was removed in~\cite{campi1998adaptive,prandini2000adaptive}.
Forced exploration techniques,~somewhat similar in spirit to the $\epsilon_t$-greedy learning algorithm of \cite{sutton2018reinforcement,auer2002finite}, are studied in \cite{lai1982least,lai1987asymptotically,chen1987optimal} vis-a-vis ensuring that they do not suffer from the insufficient exploration.
% , i.e. the asymptotic estimates of the unknown parameters converge to their true value.~More recently, inspired by the RBMLE approach, t
An adaptive LQG control based on the UCB approach, called OFULQ, is proposed in \cite{abbasi_11} which establishes a regret of $\tilde{\mathcal{O}}(\sqrt{T})$.
% , which is dubbed as optimism in the face of uncertainty (OFU). It derives bounds on its learning regret by building upon the techniques introduced in~\cite{campi_98}. The algorithm proposed therein is ``UCB-type,'' i.e. it optimizes the function $J\ust(\te)$ over a high-probability confidence ball, in order to get an optimistic estimate of $\te$, and is shown to have a regret that is bounded as $\Tilde{\mathcal{O}}(\sqrt{T})$\footnote{$\Tilde{\mathcal{O}}$ hides any dependencies that are logarithmic in $T$.}.~
A similar algorithm to address the adaptive LQG control problem is also designed in \cite{ibrahimi2012efficient}.
A computationally efficient algorithm called ROBUST with a regret of $\tilde{\mathcal{O}}(T^{\frac{2}{3}})$
is proposed in \cite{dean_17}.
%\prk{I have deleted the reference to the work of Dean that just studies sample complexity since it is not relevant here. On the other hand should we not refer to the other paper of Dean that proposes the Robust algorithm?}
% More recently,~\cite{dean2020sample} studies a related ``sample complexity" problem associated with LQG control.
% , i.e. it quantifies the number of samples required in order to generate a sufficiently precise estimate of the parameters $(A,B)$. 
% However, it does this under the assumption that the controls are generated in an i.i.d. manner (which is somewhat equivalent to forced exploration).~
An alternative approach for designing learning algorithms is Thompson sampling~\cite{thompson1933likelihood},
% , or posterior sampling. 
and \cite{abbasi2015bayesian,ouyang2017control} have established an expected regret of $\Tilde{\mathcal{O}}(\sqrt{T})$ in a Bayesian context. 
More recently, new DE-based algorithms, IP \cite{faradonbeh2018input} and RCE \cite{faradonbeh_rce}, have been proposed that make an additional assumption, not made in ARBMLE, OFULQ or TS, that one has access to a stabilizing controller for the unknown system, and
% Both fall in the broad category of DE approaches, using perturbations to the control inputs that are suggested by the certainty equivalence controller, and 
establish $\tilde{O}(\sqrt{T})$ regret.

\subsection{The RBMLE Approach}
The RBMLE approach, proposed four decades ago \cite{kumar1982new}, 
was the first approach not resorting to forced choices as in \cite{robbins1952some}.
% It suggested adding a term preferring parameters $\te$ with lower costs to the likelihood function. 
% We will describe RBMLE in a simple Markov Decision Process (MDP) setup that can be generalized to other Reinforcement Learning problems such as the LQG problem as well as other Controlled Markov Processes.
% setup that is more general than the adaptive LQG problem. This setup can be used to model reinforcement learning (RL) problems and adaptive control of nonlinear stochastic control systems.
We begin by giving an informal description of it in the context of the adaptive LQ problem. 
Since the true parameter $\te\ust := (A\ust,B\ust)^T$ is not known, one can make a Least-Squares Estimate (LSE) 
$\hat{\te}_t = (\hat{A}_t, \hat{B}_t)^{\intercal}$ of them:
\alig{ 
\hat{\te}_t \in \mbox{ArgMin}_{\te=(A,B)^{\intercal}}  \sum_{s=0}^{t-1} ||  x_{s+1}-Ax_s-Bu_s ||^2. \label{LSE}
}
Under the certainty-equivalence approach, the control input applied is $u_t = K(\hat{\te}_t)x_t$ where $K(\te)$ is the optimal linear feedback gain for the LQ problem when the system is described by $\te=(A,B)^{\intercal}$.
Suppose now that these estimates were to converge as $t \to \infty$ to 
$\hat{\te}_\infty = (\hat{A}_\infty, \hat{B}_\infty)^{\intercal}$.
Then, asymptotically, the input applied is $u_t \approx K(\hat{\te}_\infty) x_t$.
The closed-loop system therefore settles down to behaving according to
\alig{
x_{t+1} \approx (A\ust+B\ust K(\hat{\te}_\infty)) x_t + w_{t+1}.
} \label{cloop}
As it does so, one loses the
ability to identify the matrices $(A\ust,B\ust)$, and
asymptotically one can only identify the closed-loop gain $A\ust+B\ust K(\hat{\te}_\infty)$.
The problem is that as the parameter estimates begin to converge to $\hat{\te}_{\infty}$, the control gain converges to $K(\hat{\te}_\infty)$
and further exploration ceases, and one ends up only identifying the behavior of the system under the limiting gain $K(\hat{\te}_\infty)$ being applied to the system. 
% One ceases to explore for other potentially better policies. 
Since the limiting policy $u_t=K(\hat{\te}_\infty)x_t$ need not be optimal for the long-term average cost for the true system $(A\ust,B\ust)$, the CE rule leads to a sub-optimal performance. Indeed, this problem goes by various names in different fields -- the dual control problem \cite{feldbaum1960dual,wittenmark1995adaptive}, the closed-identifiability problem \cite{borkar1979adaptive,berry1985bandit,gittins2011multi}, or exploration vs. exploitation dilemma \cite{lattimore2020bandit}.

This problem was resolved in \cite{kumar1982new}
without resorting to forced exploration as in~\cite{robbins1952some}.
% \rs{Do we want to mention that the first work to do so, since Abbasi says they and Lai Robbins are the first works to do so?}~
The key observation made there was that under CE the limiting parameter estimate $\hat{\te}(\infty)$ has a one-sided bias. Specifically, $J\ust(\hat{\te}_{\infty}) \geq J\ust (\te\ust)$, i.e., the limiting parameter estimate has an optimal cost that is larger than the optimal cost of $\te\ust$.
To see this, denote by $J(K, \te)$ the long-term cost of using the control gain $K$ when the parameter is $\te$. Then the fact that the models $\hat{\te}_{\infty}$ and $\te\ust$ have identical behavior under the gain $K(\hat{\te}_{\infty})$
implies that
$J(K(\hat{\te}_{\infty}),\hat{\te}_{\infty}) = J(K(\hat{\te}_{\infty}),\hat{\te}\ust)$.
Now note that $J(K(\hat{\te}_{\infty}),\hat{\te}_{\infty})=J\ust(\hat{\te}_{\infty})$ since the gain $K(\hat{\te}(\infty))$ is optimal for $\hat{\te}_{\infty}$. However the gain $K(\hat{\te}_{\infty})$
is not necessarily optimal for $\te\ust$, and so $J(K(\hat{\te}_{\infty}),\hat{\te}\ust) \geq J\ust(\te\ust)$.
Therefore,
\begin{equation} \label{bias}
  J\ust(\hat{\te}_{\infty})  =  J(K(\hat{\te}_{\infty}),\hat{\te}_{\infty}) = J(K(\hat{\te}_{\infty}),\hat{\te}\ust) \geq J\ust(\te\ust).
\end{equation}

Following from this observation, it was reasoned in \cite{kumar1982new} that
if one could slightly bias the MLE to favor models with lower optimal costs 
$J\ust(\te)$
so  as to obtain 
% \begin{equation} \label{the-fix}
$J\ust(\hat{\te}_{\infty}) \leq J\ust(\te\ust)$,
% \end{equation}
then one would have equality throughout (\ref{bias}), to obtain
% \begin{equation} \label{opt}
 $J(K(\hat{\te}_{\infty}),\hat{\te}\ust) = J\ust(\te\ust)$,
%   J\ust(\hat{\te}_{\infty})  =  J(K(\hat{\te}_{\infty}),\hat{\te}_{\infty}) = J(K(\hat{\te}_{\infty}),\hat{\te\ust}) = J\ust(\te\ust),
% \end{equation}
yielding the desired result that the gain $K(\hat{\te}_{\infty})$ is optimal for $\te\ust$.
So motivated, \cite{kumar1982new} proposed RBMLE\footnote{It has been called the Cost-Biased MLE, as in \cite{campi1998adaptive}, or  Reward-Biased MLE (RBMLE) as in \cite{liu2020exploration}, depending on whether one minimizes a cost or maximizes a reward.} which in the context of the LSE
 suggests choosing the parameter estimate
$\hat{\te}_{RBMLE}(t) =(\hat{A}_{RBMLE}, \hat{B}_{RBMLE})^{\intercal}$ as:
\alig{ 
\hat{\te}_{RBMLE}(t) \in
\mbox{ArgMin}_{\te=(A,B)^{\intercal}} \Big[ \alpha(t) J\ust(\te)
  + \sum_{s=0}^{t-1} || x_{s+1}-Ax_s-Bu_s ||^2 \Big]. \label{BMLE}
}
Generally $\alpha(t) = o(t)$ with the precise growth rate of $\alpha(t) \nearrow +\infty$ dependent on the context.
% appropriate for an LQ system with i.i.d.~Gaussian noise, when the true parameter belongs to a known compact set in \cite{kumar1982optimal,campi1998adaptive}. 
% Specifically, it was shown that the long-term average reward is optimal, which means that the ``regret" is $o(T)$.

% This must agree with the estimated closed-loop gain
% $\hat{A}_\infty+\hat{B}_\infty K(\hat{\te}_\infty)$.
% Therefore, asymptotically, (\ref{cloop}) behaves similarly to
% \alig{
% x_{t+1} \approx (\hat{A}_\infty+\hat{B}_\infty K(\hat{\te}_\infty)) x_t + w_{t+1}.
% } \label{cloop2}
% Therefore the quadratic cost incurred by (\ref{cloop})
% must equal the quadratic cost incurred by (\ref{cloop2}).
% Now the gain $K(\hat{\te}_\infty))$
% is optimal for $(\hat{A}_\infty,\hat{B}_\infty)$ but is generally only suboptimal 
% for $(A\ust,B\ust)$. Hence if $J(\pi, \te)$ denotes the long-term average cost of using a policy $\pi$ on the model $\te$, and $J\ust(\te)$ the optimal cost achievable for the system modeled by $\te$, then one has

% optimal cost achievable for a system $\te$, one has
% \alig{
% J(\hat{\te}_\infty) \geq J(\te\ust).
% }
% Therefore there is a one-sided bias in the parameter estimates; they tend to converge to parameters that have larger optimal costs.
% This suggests that if one could delicately values 
% one can determine the cost of the using the feedback gain $K(\hat{\te}_\infty$ on the 

% In this work we analyze RBMLE algorithm in the context of RL for LQG control systems.

\subsection{The UCB Approach}
The ``Upper Confidence Bound'' (UCB) approach was first proposed in the context of Multi-Armed Bandit (MAB) Problems
in \cite{lai1985asymptotically}.
% While the precise policy designed in \cite{lai1985asymptotically} has some additional complexity,
In the context of Bernoulli bandits, it essentially consists of constructing, for each arm $i$ at each time $t$, a confidence interval $( \te(i,1;t), \te(i,2;t) )$ of its payoff probability with confidence $(1-\delta(t))$, for $\delta (t) = o(\frac{1}{t})$,
and playing the arm $i$ with the highest value of $\te(i,1;t)$.
 This approach has been 
generalized to a variety of contexts including linear contextual bandits \cite{chu2011contextual}, Gaussian Processes \cite{srinivas2009gaussian},  MDPs \cite{auer2007logarithmic,jaksch2010near}, and LQ systems \cite{abbasi_11}.

Particularized to LQ systems, the OFULQ algorithm of \cite{abbasi_11} that is based on the UCB approach suggests
choosing $\hat{\te}_{OFLUQ}(t) = (\hat{A}_{OFULQ},\hat{B}_{OFULQ})^{\intercal}$ as the minimizer of:
\alig{ 
  \min_{\te \in {\mathcal{C}_t(\delta)}} J\ust(\te), \mbox{  where} \label{UCB-objective}
}
% where
\begin{equation} \label{UCB-constraint}
 \mathcal{C}_t(\delta):=
 \left\{\theta = (A,B)^{\intercal}:
 \sum_{s=0}^{t-1} || x_{s+1}-Ax_s-Bu_s ||^2 \leq \gamma_T(\delta) \right\}
% \left\{\theta: trace\left((\theta-\hat{\theta}_t)^TZ_t(\theta-\hat{\theta}_t)\right) \leq \beta_t(\delta)\right\},\notag
\end{equation}
is a  $(1-\delta)$-high confidence set of parameters for an appropriate choice of $\gamma_T(\delta)$ \cite{abbasi_11}.

%Why then do we not simply use Least Squares Estimate, i.e. the Certainty Equivalence controls? The current state of the field is that there is no regret bound proof for unmodified Certain Equivalence based on $\theta_{LSE}$ (CE-LSE), or a counterexample that CE-LSE is unstable or has bad regret.~\cite{simchowitz20} assumes knowledge of a stabilizing gain matrix and switches to it occasionally, and so does not provide an answer this question. The ODE method studied in~\cite{lin1985will} suggests convergence to an equilibrium manifold, but does not provide information on a possible slower drift along the manifold to the true parameter.}

\section{The Augmented RBMLE-UCB Method}

The UCB approach has been called ``Optimism in the Face of Uncertainty'' (OFU), since it chooses an optimistic arm after calculating the confidence intervals. The RBMLE version is also the same,
though it does it in a different way by directly giving preference to parameters that can yield better rewards.
It arrives at optimism in a very systematic way by noticing that closed-loop identification 
% (\ref{closed-loop identification}) 
leads to the chain of inequalities (\ref{bias}), which could then be made into equalities by
ensuring that $J\ust(\hat{\te}_{\infty}) \leq J\ust(\te\ust)$.
% This is a more precise reason than just ``optimism" for preferring $\te$'s with smaller costs.

We now show how one may reconcile the two approaches to optimism in the face of uncertainty in the context of the LQ problem, and then combine them to obtain a method that has experimentally superior empirical performance while also allowing a proof that it achieves the currently best known order of regret.
% Based on that we design a develop that combines the two approaches coherently.
 
% We begin by noticing that particularized to the context of the LQG problem, RBMLE suggests choosing the parameter estimate 
% $\hat{\te}_{RBMLE}(t) =(\hat{A}_{RBMLE}, \hat{B}_{RBMLE})$ as the minimizer of:
% \alig{ 
% \min_{\te} \Big[ \alpha(t) J\ust(\te)
%   + \sum_{t=0}^{T-1} \left(  x_{t+1}-Ax_t-Bu_t \right)^2 \Big], \label{BMLE}
% }
% since the log-likelihood is the sum of squares given the Gaussian distribution
% of the noise (though we also analyze the sub-Gaussian case in the sequel).

It can be seen that the RBMLE approach (\ref{BMLE}) can be considered an unconstrained penalty version of the constrained optimization problem \eqref{UCB-objective}-\eqref{UCB-constraint}, with a penalty factor $\frac{1}{\alpha(t)}$ for constraint violation.
This provides a justification for employing ``optimism'' in the UCB approach.
% , as another way of
% converting (\ref{bias}) into a chain of equalities.

This raises the question of whether we can take advantage of these synergies to fashion a superior algorithm,~that is superior from the viewpoint of being able to provide a theoretical guarantee of the best known order of regret to date, as well as superior from the point of view of providing the best experimental performance of the algorithms to date. One can draw inspiration here from the Augmented Lagrangian Method of \cite{hestenes1969multiplier,powell1969method}. Given a constrained problem as in 
\eqref{UCB-objective}-\eqref{UCB-constraint}, it adds a penalty for constraint violation, and also adjoins the constraint through a Lagrange multiplier, thus involving the constraint twice.
% This methods performs very efficiently and can even obtain super-linear convergence over macro steps of the dual update. 
In our case, we add the penalty and also retain the constraint.
% one can just treat the problem in the primal space by not adjoining the constraints.
This leads to the Augmented RBMLE-UCB method:
\alig{ \hat{\te}_{\mbox{Aug}}(t) \in
\mbox{ArgMin}_{\te \in {\mathcal{C}_t(\delta)}} \Big[ \alpha(t) J\ust(\te) +
   \sum_{s=0}^{t-1} || x_{s+1}-Ax_s-Bu_s ||^2 \Big]. 
}
% \alig{
% \min_{\te \in {\mathcal{C}_t(\delta)}} \Big[ \alpha(t) J\ust(\te)
%   + \sum_{t=0}^{T-1} \left( x_{t+1}-Ax_t-Bu_t \right)^2 \Big].
% }
% This algorithm has the following behavior. As long as the minimizer is on the boundary of $C_t(\delta)$, the estimate coincides with what UCB would choose. However, if the minimizer is in the interior of $C_t(\delta)$, then it chooses what $\hat{\te}_{\mbox{RBMLE}}(t)$ would locally choose.
% Another 
The advantage of retaining the constrained optimization is in enabling theoretical analysis, in that we can make use of the
bounds on parameter estimates from concentration inequalities, as we shall 
show in the sequel.
We thereby prove in Section \ref{sec:regret}
that it also has $\Tilde{\mathcal{O}}(\sqrt{T})$ order of regret, which is the best known so far \cite{abbasi_11}. Moreover, simulations reported in Section \ref{sec:simulations} show that its performance is also much improved. In all cases, the simulation performance of ARBMLE is the same as RBMLE, raising the open problem of proving that the original RBMLE for the LQ problem also has the best known order of regret.

\section{Problem Formulation}
\label{sec:formulation}
As introduced in section 
\ref{sec:intro}, we consider the following linear system,
\begin{equation}
x_{t+1}=A\ust x_t+B\ust u_t+w_{t+1}
\end{equation}
where, $A\ust \in \mathds{R}^{n \times n}$ and $B\ust \in \mathds{R}^{n \times m}$ are unknown system parameters. Define $z^\top _t:= (x_t^\top , u_t^\top )$.
% \begin{pmatrix}x_t \\ u_t\end{pmatrix}$ 
Then, the linear stochastic system can equivalently be written as 
\begin{equation}
    x_{t+1}={\theta\ust}^\top z_t+w_{t+1}.
\end{equation}
The system incurs a cost at time $t$ given by,~$x^{\intercal}_t Qx_t+u^{\intercal}_t Ru_t$.~We assume that $Q \in \mathds{R}^{n \times n}$ and $R \in \mathds{R}^{m\times m}$ are known positive semi-definite and positive definite matrices respectively.
We make the following assumptions on $\theta\ust= (A\ust, B\ust)^{\intercal}$.
\begin{assumption}
\label{assumption:parameters}
There exists a known positive constant $c$
such that ${\theta\ust} \in \mathcal{S}=\mathcal{S}_0 \cap \mathcal{S}_1$ where,
% \begin{align*}
        $\mathcal{S}_0 =\{ \theta=\left(A,B\right)^{\intercal} | \left(A,B\right) \text{ is stabilizable and} \left(A,Q^{1\slash 2}\right) \text{ is detectable} \}$\footnote{We note that even though~\cite{abbasi_11} requires each $\theta$ from $\mathcal{S}_0$ to satisfy the stronger reachability-observability assumption, the results therein hold true if this is replaced by the weaker stabilizability-detectability assumption~\cite[p.61]{abbasi2013online_thesis}. In fact all that is required for the proofs to go through, is that for each $\te\in \mathcal{S}_0$, the corresponding Riccati equation %\prk{\sout{to}}
        % ~\eqref{eq:riccati_1} 
        have a unique solution. It is well-known~\cite{lewis2012optimal} that stabilizability and detectability, as above, is sufficient for this.} \cite{kailath1980linear},
     and
        ${\cal{S}}_1 = \{ \theta=\left(A,B\right)^{\intercal} | trace(\theta^{\top}\theta)\leq c \}$.
% \end{align*}
\end{assumption}

\begin{assumption}\label{assumption:noise}
Let $\mathcal{F}_t :=\sigma((x_0,u_0),(x_1,u_1),\cdots,(x_t,u_t))$
denote the history of states and inputs until time $t$.
The state noise $\{w_t\}$ is a martingale difference sequence with respect to $\{\mathcal{F}_t\}$, with
$\mathds{E}[w_{t+1}w^{\intercal}_{t+1}|\mathcal{F}_t]=I_n$, and is element-wise sub-Gaussian, i.e., for any $\gamma \in \mathds{R}$, there exists $L>0$ with
% \begin{equation}
    $\mathds{E}\left[exp(\gamma w_{t,j})|\mathcal{F}_t\right] \leq exp\left(\frac{\gamma^2L^2}{2}\right)~ \forall ~j \text{ and } \forall~ t$.
% \end{equation}

% Consider the discrete time linear system~\cite{kumar_15},
% \begin{equation}\label{eq:system}
%     x_{t+1}=A\ust x_t+A\ust u_t+w_{t+1}
% \end{equation}
% where $x_t \in \mathds{R}^n,u_t \in \mathds{R}^m$ are the state and the control input respectively at time $t$. $w_t \in \mathds{R}^n$ represents the noise in the system at time $t$. The system parameters, represented by matrices $A\ust \in \mathds{R}^{n \times n}$ and $B\ust \in \mathds{R}^{n \times m}$ are unknown.

%We denote \prk{If you don't ever use $c_t$ anywhere in the paper, then you can delete this sentence and equation.} 

%\prk{$R$ should be a positive definite matrix, not enough to be semidefinite}.
% while for a vector $x$, we let $x^\top $ denote its transpose.
% The sample-path average cost incurred by the controller that chooses $u_t$ at time $t$ is: %\prk{Probably better to make lim below into limsup since limit may not exist}
% %$$
% %   J(x(0),u(0), \dots) = 
% %$$
% \begin{equation} \label{eq:avg_cost}
% \limsup_{T \to \infty}\frac{1}{T} \sum_{t=0}^\top  c_t.
% \end{equation}
% \subsection*{Notations:} 
% Let $I_n$ be identity matrix of dimension $n \times n$ for any $n \in \mathds{N}$.
\end{assumption}
For a matrix $M$, we use $\|M\|$ to denote its operator norm induced from the $\ell_2$-norm.
\subsection{Controller Design for a Known LQG System}
For a discrete time linear system,
$x_{t+1}=\theta^{\intercal} z_t+w_{t+1},$
where $\theta = (A,B)^{\intercal}  \in \mathcal{S}$, there is an unique positive semidefinite matrix $P(\te)$ %\footnote{We suppress its dependence upon the parameter $\te$})
that satisfies the Riccati equation (see \cite{kumar1986stochastic})
% \begin{equation}\label{eq:riccati_1}
   $$P(\te)= Q+A^{\intercal} P(\te)A-A^{\intercal} P(\te)B(B^{\intercal} P(\te) B+R)^{-1}B^{\intercal} P(\te)A .$$
% \end{equation}
The optimal control law which minimizes the long-term average quadratic cost
% ~\eqref{eq:avg_cost} 
is $u_t=K(\theta)x_t$, where the ``gain matrix'' is 
% \begin{equation}
$K(\theta):=-(B^{\intercal} P(\te)B+R)^{-1}B^{\intercal} P(\te)A$.
% \end{equation} 
The optimal average cost, $J\ust(\theta)$ is equal to $trace(P(\theta))$. 
As a consequence of Assumption \ref{assumption:parameters}, the parameter set $\mathcal{S}$ is bounded, hence one can show that $P(\theta)$ is bounded as well: 
\begin{equation}
        D:= \sup_{\theta \in \mathcal{S}} ||P(\theta)||  < \infty.\label{def:D}
    \end{equation}
The following assumption is commonly made in online LQR learning problem in which the knowledge of a stabilizing controller is not made~\cite{abbasi_11,ouyang2017control}.
\begin{assumption}\label{assum:bounds}
    \begin{equation}
    \rho:= \sup_{\theta \in \mathcal{S}} ||A+BK(\theta)|| <1, \mbox{  and}\label{eq:assum_1}
    \end{equation}
    % There exists a positive constant
    % $c_0$ such that, 
    \begin{equation}
        c_0:= \sup_{\theta \in \mathcal{S}} ||K(\theta)|| < +\infty.
    \end{equation} 
% \end{enumerate}
\end{assumption}
Note that there are recent works which relax~\eqref{eq:assum_1} slightly to $\sup_{\theta \in \mathcal{S}} \rho(A+BK(\theta))<1$, however they assume that the learner has access to a stabilizing controller \cite{dean_17,dean_20,simchowitz20}.
%Under assumption \ref{assumption:parameters}
%, the following results hold true:
%\begin{enumerate}
%    \item Matrix $A+BK(\theta)$ is stable, i.e.,  $||A+BK(\theta)||_2 < 1$. \prk{Do you rather mean all eigenvalues inside the unit circle?} \am{Yes} \prk{But saying $\ell_2$ norm is less than one is not the same as saying that all the eigenvalues are less than one.}  Let 

%    \item There exist a finite constant $D$ such that $D=\sup_{\theta \in \mathcal{S}} \left\{||P(\theta)|| \right\} < \infty.$
%\end{enumerate}
%\end{lemma}

\subsection{Construction of Confidence Interval}
\label{sec:lse}

%The squared fitting error of $\theta$ is defined as $
% \begin{equation}
%   V_t(\theta):=  $.
% \end{equation}
The $\ell_2$-regularized squared fitting error with parameter $\lambda>0$ is given by:
 \begin{equation}
   V_t(\theta)=\lambda||\theta||_2^2+\sum_{s=0}^{t-1} ||x_{s+1}-\theta^{\intercal} z_s||_2^2.
 \end{equation}
Let $\hat{\theta}_t$ be the $\ell_2$ regularized least-squares estimate of $\theta\ust$, i.e., $
% \begin{equation}
\hat{\theta}_t \in \arg \min_{\theta \in \mathcal{S}} V_t (\theta)$.
% \end{equation}
Next,
given the history $\mathcal{F}_t$, we construct a ``high-probability confidence ball,'' $\mathcal{C}_t(\delta)$ around $\hat{\theta}_t$, i.e., a set of plausible system parameters that contains the true parameter $\theta\ust$ with a high probability. 
% Let $\mathcal{C}_t(\delta)$ be the following confidence ball around $\hat{\theta}_t$,
Let 
 \begin{equation}
 \mathcal{C}_t(\delta):=
\left\{\theta: trace\left((\theta-\hat{\theta}_t)^{\intercal} Z_t(\theta-\hat{\theta}_t)\right) \leq \beta_t(\delta)\right\},
 \end{equation}
where $Z_t :=\lambda I_{n+m} +\sum_{s=0}^t{z_sz_s^{\intercal}}$, and
\begin{equation}\label{eq:conf}
    \beta_t(\delta):=\left(nL\sqrt{2\log \left(\frac{\sqrt{det(Z_t)det(\lambda I)}}{\delta}\right)}+\sqrt{\lambda}c
    \right)^2.
\end{equation}

Define 
\begin{align}
\cE_1(t)= \left\{\theta\ust \in \mathcal{C}_s(\delta/4), ~\forall~ s=1,2,\cdots,t\right\} \text{and }
\cE_2(t)= \left\{||x_s||\leq d_t,~ \forall s \leq t\right\}, \label{E2}
\end{align}
where $d_t$ is defined in the Appendix. Let $\cE_1:=\cE_1(T)$ and $\cE_2:=\cE_2(T)$.~Then, as in
\cite{abbasi_11}, $
% \begin{equation} \label{lemma:event}
    \mathds{P}\left(\cE_1 \cap \cE_2 \right) \geq 1 -\frac{\delta}{2}$.
% \end{equation}
% \end{lemma}
% Lemma 4 in \cite{abbasi_11}.
Moreover, on
% \begin{lemma}\label{lemma:x}
the event $\cE_2(t)$ defined in (\ref{E2}), the following holds, 
%\prk{I guess the point is that it is even less than $d_t$. So is the basic idea here that by  assuming a gross bound you can prove a better bound? Also, wat are $c_1$ and $c_2$?}
\begin{align}\label{eq:x}
\max_{1 \leq s \leq t} \| x_s\|\leq X_t =Y_t^{n+m+1},~\forall t=1,2,\ldots, \text{where,}
\end{align}
%\prk{Should this rather be $max_{1 \leq s \leq t} \| x_s\| =: X_t \leq Y_t^{n+m+1}$?}
\begin{align*}
    & Y_t:= \max\bigg(e,\lambda(n+d)(e-1) \times \\
    & 4\left(c_1\log\left(\frac{1}{\delta}\right)+c_2\log\left(\frac{t}{\delta}\right)\right)\bigg)\times \left( \log^2\left(4c_1\log \left(\frac{1}{\delta}\right)+4c_2\log \left(\frac{t}{\delta}\right)\right) \right)  ,
\end{align*}
with $c_1$ and $c_2$ being problem dependent constants.
% \begin{proof}
% Lemma 5 in \cite{abbasi_11}.
% \end{proof}
% \end{lemma}

\section{The Augmented RBMLE-UCB Algorithm}
\label{sec:algo}
We employ a version of the Augmented RBMLE-UCB (ARBMLE) algorithm that proceeds in an episodic manner. Let $t_k$ denote the starting time of the $k$-th episode. Then, during episode $k$, it implements the control policy,
$u_t=K(\theta_{t_k})x_t, ~\forall ~t \in \left\{t_k,t_{k}+1,\ldots,t_{k+1}-1\right\}$.
where, $\theta_{t_k}$ is obtained by solving the following optimization problem,   
\begin{equation}\label{eq:rbmle}
    \theta_{t_k} \in \arg \min_{\theta \in \mathcal{S} \cap \mathcal{C}_{t_k}(\delta)} \left\{V_{t_k}(\theta)+ \alpha(t_k)J^{\star}(\theta)\right\},
\end{equation}
where the bias-term, $\alpha(t)=\alpha_0\sqrt{T},\forall t$, for $\alpha_0 > 0$.
\begin{algorithm}[h]
   \caption{Augmented RBMLE-UCB (ARBMLE)}
   \label{algo}
\begin{algorithmic}
   \STATE {\bfseries Initialize:} $t=0$,~$Z_0=\lambda I_{n+m}$
      \FOR{$k=0,1,\cdots$}
    \IF{$det(Z_t)>2det(Z_{t_{k-1}})$} 
  \STATE solve the following optimization to obtain $\theta_{t_k}$,
 \STATE \begin{equation*}
      \theta_{t} \in\arg\min_{\theta \in \mathcal{S} \cap \mathcal{C}_{t_k}(\delta)} \left\{V_{t_k}(\theta)+ \alpha(t_k)J\ust(\theta)\right\},
  \end{equation*}
    \ELSE 
    \STATE $\theta_t=\theta_{t-1}$
     \ENDIF
        \STATE $u_t=K(\theta_t)x_t$
        \STATE $Z_{t+1}=Z_t+z_t z^{\intercal}_t$
        \STATE $t \to t+1$
   \ENDFOR
\end{algorithmic}
\end{algorithm}

%\section{Regret Analysis }
\label{sec:regret}
% We focus on the theoretical analysis of the Augmented RBMLE-UCB algorithm defined in Section \ref{sec:algo}. More specifically, we derive an upper bound on the regret of the algorithm:
% % as defined as~\cite{abbasi_11}: 
% \begin{align}
%   \label{def:regret}
%     R(T)=\sum_{t=0}^T c_t - TJ\ust(\theta\ust).
% \end{align}
% The notion of regret was first introduced in \cite{lai1985asymptotically} in context of multi-armed bandit problem. It has been widely used for characterizing the finite-time performance of algorithms that make sequential decisions~\cite{lattimore2020bandit,auer2002finite}. Note that minimizing the regret is equivalent to minimizing the performance cost.
%In this section, we provide an upper bound on the regret $R(T)$. In Lemma \ref{lemma:decomposition},
%we decompose the regret $R(T)$ into easy-to-analyze terms. Then, we provide upper bounds on each of these four terms, which allows us to get a upper bound on the regret.

% \begin{proof}
% See Appendix \ref{apndix:decom}.
% \end{proof}
%We bound $R_1,R_2,R_3$ and $R_4$ on the event $E_1 \cap E_2$. 
% \subsection{Bound on $R_1$}
%\prk{It would be better to call them as Lemmas rather than Theorems. Please correct future references so that they correspond to the following Lemma numbers.}

% \begin{proof}
% See Appendix \ref{apndix:r4}.
% \end{proof}
% Now that we have established upper bounds on the terms $R_1,R_2,R_3$ and $R_4$, we can use these 
%We now use these bounds in conjunction with 
% the decomposition result of 
%Lemma~\ref{lemma:decomposition} to bound the regret of Augmented RBMLE-UCB.
In Theorem \ref{thm:r}, we show that regret for ARBMLE is upper bounded by $\Tilde{\mathcal{O}}\left(\sqrt{T\log \frac{1}{\delta}}\right)$ which is same order as OFULQ \cite{abbasi_11}. 
\begin{theorem}\label{thm:r}
For any $\delta \in (0,1)$ and $T>0$, with a probability at least $(1 - \delta )$, the regret of the ARBMLE Algorithm is upper-bounded by $R(T)\leq  \Tilde{\mathcal{O}}\left(\sqrt{T\log \frac{1}{\delta}}\right)
$.
\begin{proof}
Appendix \ref{apndix:r}.
\end{proof}
%\prk{It is not good form to write $a\leq O(b)$, because $a = O(b)$ means the same thing, and in fact the $O$ notation was defined precisely for this purpose. Actually, in this case if pre-constants are important, it may be better to say $R(T)\leq c \sqrt{T} $, and say what $c$ is. }
% \begin{equation}
%         R(T)= \Tilde{\mathcal{O}}\left(\sqrt{T\log \frac{1}{\delta}}\right).
% \end{equation}
% where problem dependant constants and terms logarithmic in $T$ are hidden. 
%\prk{I changed it to $\Tilde{\mathcal{O}}$ instead of $O$ and deleted the mention of logarithmic terms.}
\end{theorem}

% \begin{proof}
% See Appendix \ref{apndix:r}.
% \end{proof}
\section{Empirical Performance}
\label{sec:simulations}
We evaluate the empirical performance of ARBMLE as well as
standard (unaugmented) RBMLE. We compare these algorithms with
OFULQ \cite{abbasi_11}, Thompson Sampling (TS) \cite{abeille_17}, Input Perturbations (IE) \cite{faradonbeh2018input}, Randomized Certainty Equivalence (RCE) \cite{faradonbeh_rce}, and Stabl \cite{lale22}.
The results shown here are for the following examples of linear systems that have appeared in the recent literature on adaptive control of linear systems: \begin{enumerate}
    \item Unstable Laplacian dynamics \cite{dean_17,lale22,abbasi_19}.
    %\item Chained Integrator Dynamics \cite{dean_17}.
    \item Large transient dynamics \cite{dean_17}.
    \item Unmanned Aerial Vehicle (UAV) \cite{Zhao21,lale22}.

        \item Longitudinal Flight Control of Boeing 747 \cite{lale22}.
\end{enumerate}
The details of these examples are provided in the Appendix.

Each simulation experiment is performed for a time horizon of 500 steps,
and repeated 50 times. The reported results are the averaged values over the 100 runs.
In Figure 1, we compare the $\log _{10}$ of the averaged regret of ARBMLE, OFULQ, TS and Stabl. In Figure 2, we plot the averaged regret of RBMLE, ARBMLE, IP and RCE. We summarize the results using the regret values at $T=500$,
averaged over the 50 runs, in Table \ref{table:2}. Results for more examples are provided in the Appendix.~Details of the implementations can be found in the Appendix.

\begin{figure}[t]
     \centering
     \begin{subfigure}[b]{0.45\textwidth}
         \centering
         \includegraphics[width=\textwidth]{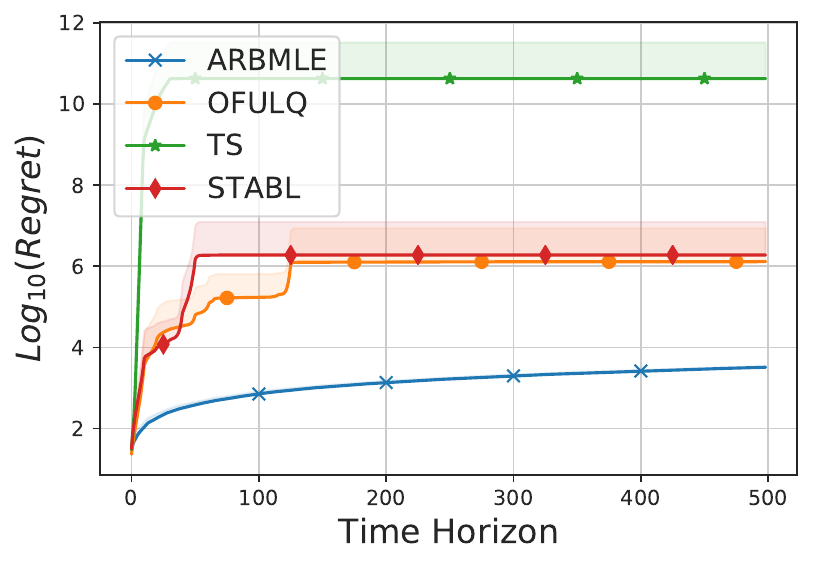}
         \caption{Unstable Laplacian}         
         \label{fig:1_3}
     \end{subfigure}
     \hfill
     \begin{subfigure}[b]{0.45\textwidth}
         \centering
         \includegraphics[width=\textwidth]{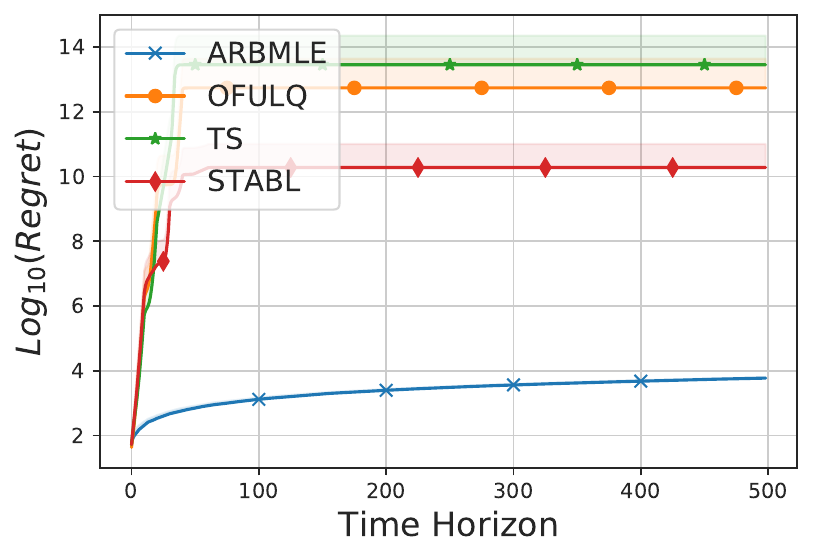}
         \caption{Large Transient Dynamics}
                  \label{fig:1_4}
     \end{subfigure}
     \hfill
          \begin{subfigure}[b]{0.45\textwidth}
         \centering
 \includegraphics[width=\textwidth]{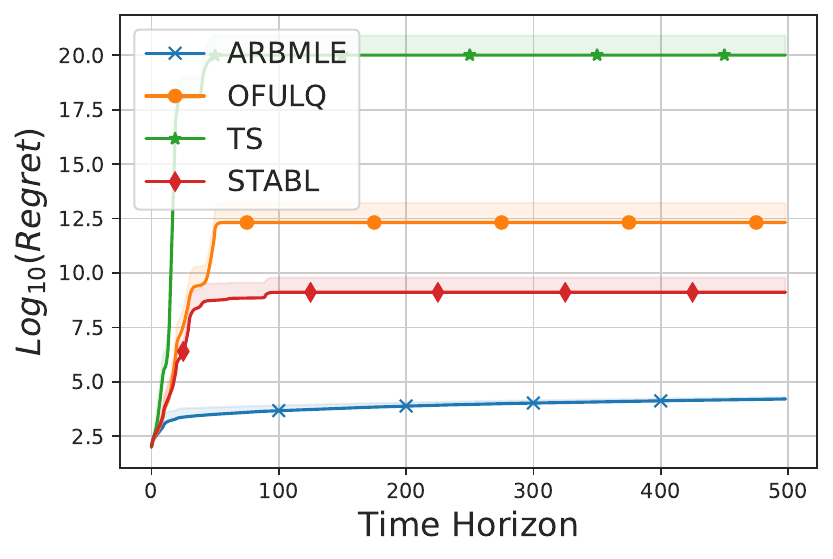}
         \caption{UAV }
              \label{fig:1_2}
    \end{subfigure}
    \hfill
 \begin{subfigure}[b]{0.45\textwidth}
         \centering
         \includegraphics[width=\textwidth]{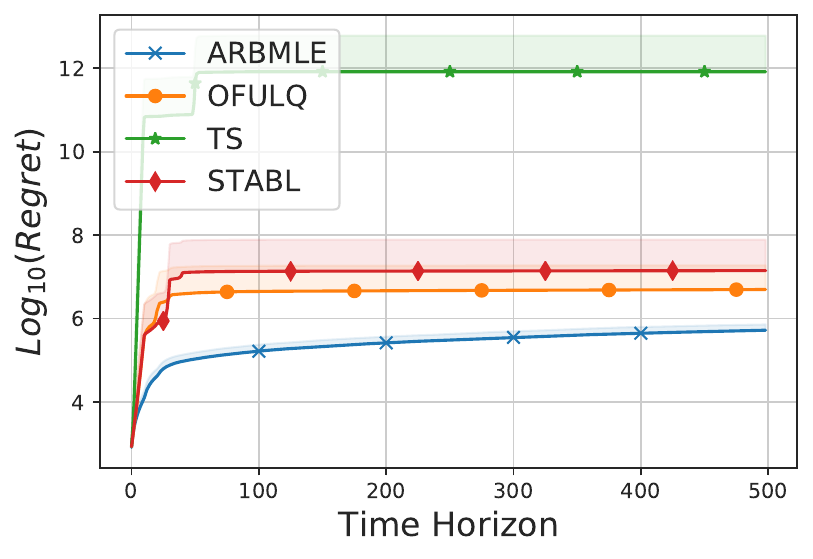}
         \caption{Boeing 747}
         \label{fig:1_1}
     \end{subfigure}
        \caption{Logarithm of the Averaged Regret over 50 runs of RBMLE, ARBMLE, OFULQ and TS for various example systems.}
        \label{fig:1}
\end{figure}

\begin{figure}[t]
     \centering
     \begin{subfigure}[b]{0.45\textwidth}
         \centering
         \includegraphics[width=\textwidth]{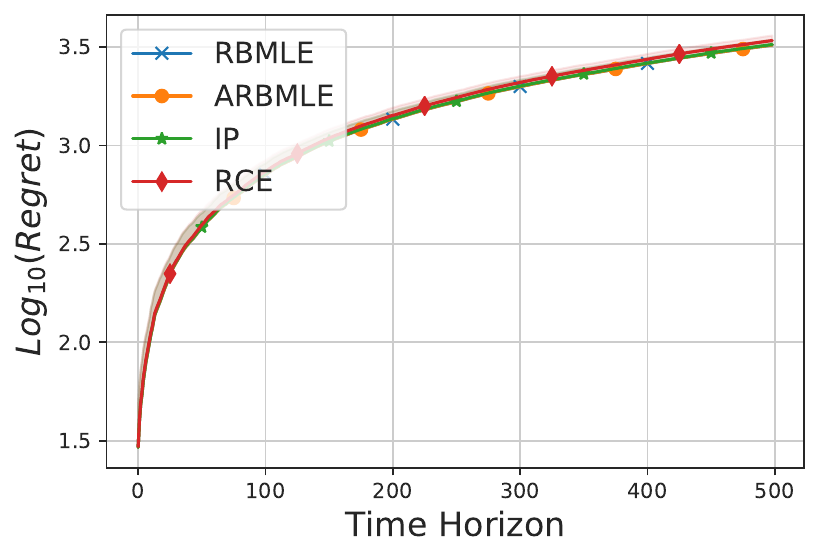}
         \caption{Unstable Laplacian}         
         \label{fig:2_3}
     \end{subfigure}
     \hfill
     \begin{subfigure}[b]{0.45\textwidth}
         \centering
         \includegraphics[width=\textwidth]{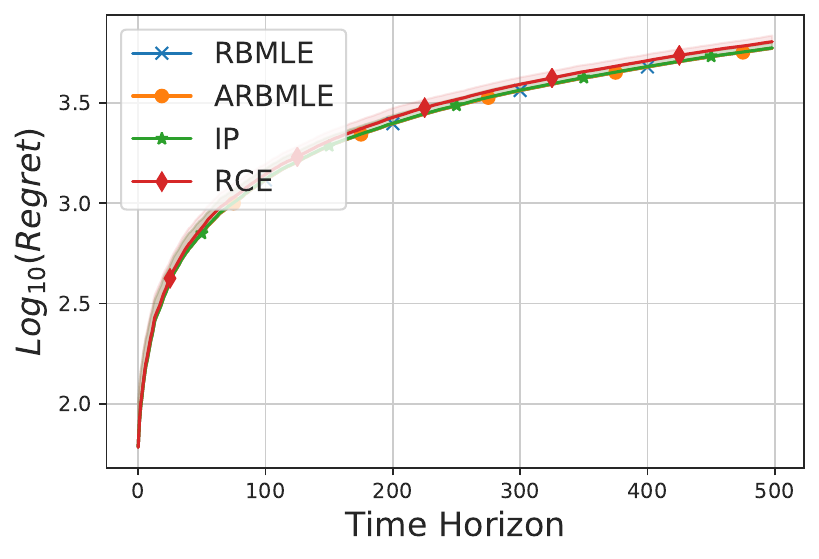}
         \caption{Large Transient Dynamics}
                  \label{fig:2_4}
     \end{subfigure}
     \hfill
          \begin{subfigure}[b]{0.45\textwidth}
         \centering
         \includegraphics[width=\textwidth]{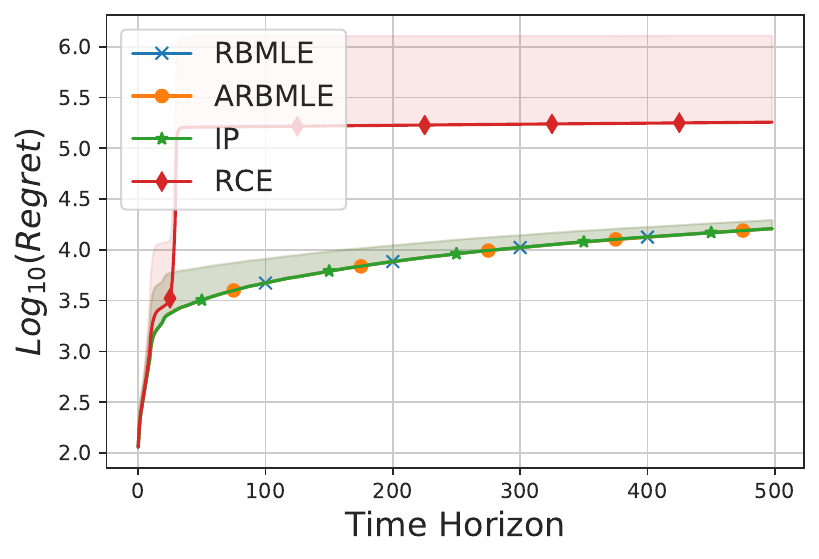}
         \caption{UAV }
         \label{fig:2_2}
     \end{subfigure}
     \hfill
          \begin{subfigure}[b]{0.45\textwidth}
         \centering
         \includegraphics[width=\textwidth]{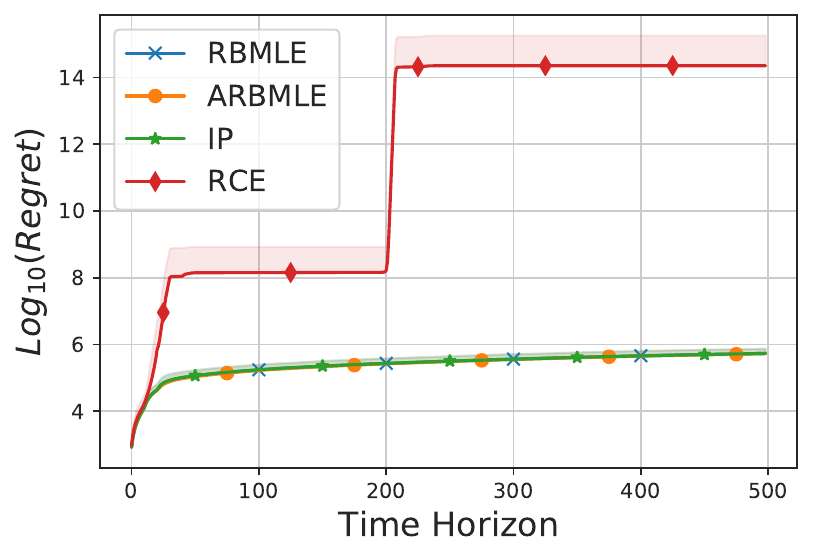}      \caption{Boeing 747}         
         \label{fig:2_1}
     \end{subfigure}
        \caption{Logarithm of the Averaged Regret over 50 runs of RBMLE, ARBMLE, IP and RCE for various example systems. Table \ref{table:2} provides the quantative results for a finer comparative evaluation of the regrets. }
        \label{fig:2}
\end{figure}

We highlight the key observations from above experiments:
\begin{itemize}[leftmargin=*]
    \item ARBMLE and RBMLE were always found to have the same empirical performance in most cases. More specifically, both these algorithms choose the same estimate $\theta_{t}$, and this suggests that $\frac{1}{\alpha(t)}$ needs to be greater than the Lagrange multiplier for the ball constraint $\theta \in \mathcal{C}_{t_k}(\delta)$.~This also motivates future study of the regret of the (unaugmented) RBMLE.
    \item As can be seen in Figure~\ref{fig:1}, ARBMLE\slash RBMLE outperform OFULQ and TS by a huge margin.~We conjecture that this is due to temporary instability when OFULQ, TS are employed~\cite{lale22}.~ARBMLE\slash RBMLE also has a signficantly lower regret as compared with StabL. 
    \item Figure~\ref{fig:2} shows that ARBMLE\slash RBMLE also outperforms RCE moderately, and IP marginally.
    % \footnote{\rs{We depict these comparisons using different plots since the regret of OFU\slash TS is much greater than that of ARBMLE\slash RBMLE\slash IP\slash RCE, and hence it is difficult to compare RBMLE with IP, RCE on the same plot.}}

\end{itemize}
\begin{center}
\begin{table}[t]
    \centering
\begin{tabular}{|c |c | c | c| c|c| c|c|} 
\hline
Ex. & RBMLE & ARBMLE & OFULQ & TS & IP & RCE & STABL\\
 \hline

(a)& 3233 & 3233 &  $1.2 \times 10^6$&  $4.2 \times 10^{10}$  & 3251 & 3408 &$1.8 \times 10^6$\\ \hline
%(c)  &2328 & 2328 & 2756 & 2580  & 2342 & 2399 &3646\\ \hline
(b)  &5930 & 5930 & $5.4 \times 10^{12}$ & $2.8 \times 10^{13}$ & 5955 & 6396 & $1.9 \times 10^{10}$\\ \hline
(c)& 16144& 16135 & $2.1 \times 10^{12}$ & $1.1 \times 10^{20}$ &16164  & 180639 & $1,2 \times 10^{9}$  \\ \hline
 (d) & 540297 & 528805 &  $4.9 \times {10^6}$&  $8.2 \times 10^{11}$  &   540248 & $2.2 \times 10^{14} $& $1.4 \times 10^{7} $ \\ \hline
\end{tabular}
    \caption{Average Regret Performance at $T=500$.}
    \label{table:2}
\end{table}
\end{center}
\begin{remark}
 Our simulations for ARBMLE, OFULQ and TS are based on the confidence interval $\beta_t(\delta) $ as defined in \eqref{eq:conf}.
%  \prk{Can you add that confidence interval is from the paper by Dean et al?}. 
 Instead, recent works \cite{lale22,dean_17} use  $\beta_t(\delta) := trace\left((\theta\ust-\hat{\theta}_t)^{\intercal} Z_t(\theta\ust-\hat{\theta}_t)\right)$.  However, one may note that $\theta\ust$ in $\beta_t(\delta)$ is not known to the learning agent, and so such a definition of $\beta_t(\delta)$ is not a viable implementation. The effect of the choice of $\beta_t(\delta)$ on the regret performance is shown in the appendix. 
 \end{remark}
\section{Concluding Remarks}
\label{sec:discussion}
We reconcile the RBMLE and UCB approaches by showing that RBMLE is an unconstrained
penalty version of the constrained optimization problem of UCB. ~Showing that UCB is a constrained version of RBMLE also explains why the optimism embodied in UCB-based schemes is justified. In particular, it is justified by the goal of nullifying the one-sided bias that
results from the closed-loop identification of dynamics.

Building on this, we propose an Augmented RBMLE-UCB method that not only matches the best known order of regret, $\tilde{\mathcal{O}}(\sqrt{T})$, to date --  that of the OFULQ algorithm -- but also outperforms OFULQ, TS and Stabl by a significant margin in simulation experiments.~In fact, for the UAV experiment, the regret of OFULQ is $\approx 10^{8}$, TS $\approx 10^{16}$ and StabL $\approx 10^5$ times that of RBMLE. It outperforms RCE also, but the gains are moderate, while with respect to IP the gains are smaller. The simulations were carried out on real-world examples such as Boeing 747, UAV, Unstable Laplacian and Large Transient Dynamics which were taken from recent works~\cite{dean_17,dean_20,abbasi_19,lale22,Zhao21} on adaptive LQG control.
% The simulations of the systems studied show that at least in these systems, the regret performance of Augmented RBMLE appears to better than the other algorithms studied.

This work further extends recent studies of RBMLE for MDPs~\cite{mete2021reward}, stochastic Multi-Armed Bandits~\cite{liu2020exploration} and Contextual bandits \cite{hung2020reward} establishing state-of-art regret results as well as empirically good performance.
% , this work on LQ systems lends more evidence that RBMLE could be a promising and competitive Reinforcement Learning algorithm. 

There remain several open questions.

Currently we study only an augmented form of RBMLE, i.e. ARBMLE, in which the agent searches for a parameter value that optimizes a certain reward-biased maximum likelihood objective function within a ``high-confidence ball.'' However, (unaugmented) RBMLE optimizes this over the larger set $\mathcal{S}$ (Assumption~\ref{assumption:parameters}) that is known to contain the true parameter. The reason for this augmentation is simply that currently we are unable to prove regret bounds for RBMLE without it. By including it, we can capitalize on the nice technical results in~\cite{abbasi_11} in the analysis of OFULQ. Simulations show that the constraint is loose, and performance is the same with\slash without the constraint. In fact, as shown by the simulations, the choices of $\theta_{t_k}$ made by the standard RBMLE, and the Augmented RBMLE are the same. It remains to be seen whether similar regret guarantees can be obtained by the (unaugmented) RBMLE, which remains an open problem. 

% Perhaps the solution of~\eqref{eq:rbmle} lies within the confidence interval. 

% Surprisingly, this also shows that a local minimum for the Augmented RBMLE is also a local minimum for standard RBMLE, even though the corresponding global minima might be different. \rs{What does this mean?:} Hence the conjecture might not be true. 

One may note that ARBMLE and OFULQ differ in their choices of decisions. Denoting their estimates of the unknown parameter at a given time by $\theta_{\text{ARBMLE}}$ and $\theta_{\text{OFULQ}}$ respectively, $\theta_{\text{OFULQ}}$
always lies at the boundary of the UCB-ball constraint.  This is easily seen to be true for bandits and MDPs, and empirically observed to be so for LQ systems. In contrast, $\theta_{\text{ARBMLE}}$ is most often in the interior of the UCB-ball.~Moreover, while OFULQ treats all models within the UCB-ball equally and only assesses them by their cost,~ARBMLE prefers models that are closer to the Least Squares Estimate ($\theta_{\text{LSE}}$). For example, if $J^{\star}(\theta_1) = J^{\star}(\theta_2)$, then ARBMLE prefers the $\theta_i$ that is closer to $\theta_{\text{LSE}}$.~We conjecture that this is the reason why ARBMLE has a significantly better performance than OFULQ -- but we are unable to prove it.

    It has been shown in previous works on adaptive LQG control~\cite{kumar1983optimal,campi1998adaptive} that under appropriate formulations RBMLE is guaranteed to stabilize an unknown LQG system over an infinite horizon, i.e., $\limsup_{T\to\infty}\frac{1}{T}\sum_{t=1}^{T} ( \|x(t)\|^{2} + \|u(t) \|^2) <\infty$ a.s..
    Morover, the sample path performance cost $\limsup_{T\to\infty}\frac{1}{T}\sum_{t=1}^{T}(x^{\intercal}(t)Qx(t) + u^{\intercal}(t)Ru(t))$ is (a.s.) equal to the optimal performance that could be attained if the system parameters $(A\ust,B\ust)$ were known. In fact, stability is a prerequisite before one can establish the latter type of result, as the decades of work on stochastic adaptive control in the nineteen seventies to the nineties has shown. While recent work has tended to study the finer performance measure of regret, it is usually over a finite horizon, and more attention to the stability of the learning process over an infinite horizon appears well deserved. Similarly, robustness which also was rigorously formulated and addressed in earlier decades needs to be re-examined \cite{praly1989robust}.%Similarly, there have also been studies 
    % by quantifying the deviations from this ``optimal performance.''  
    
\section*{Acknoweldgments
}
This material is based upon work partially supported by the
US Army Contracting Command under W911NF-22-1-0151,
US Office of Naval Research under N00014-21-1-2385;
4/21-22 DARES: Army Research Office W911NF-21-20064
US National Science Foundation under CMMI-2038625, 
The views expressed herein and conclusions contained in this document are those of the authors and should not be interpreted as representing the views or official policies, either expressed or implied, of the U.S. Army Contracting Command, ONR, ARO, NSF, or the United States Government.
The U.S. Government is authorized to reproduce and distribute reprints for Government purposes notwithstanding any copyright notation herein. The work of Rahul Singh was partially supported by the SERB Grant SRG/2021/002308, and PC 39010B.

\medskip

\bibliography{references}
\bibliographystyle{unsrt}

%%%%%%%%%%%%%%%%%%%%%%%%%%%%%%%%%%%%%%%%%%%%%%%%%%%%%%%%%%%
\clearpage
\appendix
%\section*{Supplementary Material: Augmented RBMLE-UCB Approach for Adaptive Control of Linear Quadratic Systems}
\section*{Appendix}
\pagenumbering{gobble}
\section{Regret Analysis}
We now prove the $\Tilde{\mathcal{O}}\left(\sqrt{T\log \frac{1}{\delta}}\right)$ upper-bound on regret of Augmented RBMLE-UCB that was claimed in Section~\ref{sec:regret}.  
\begin{lemma} \label{lemma:decomposition}
The regret $R(T)$ of the Augmented RBMLE-UCB learning algorithm
% which implements $u_t=K(\theta_t)x_t$
can be decomposed as $
% follow:
% \begin{equation}
    R(T)= R_1 +R_2 + R_3 + R_4$, where
    % \mbox{ where}\label{eq:regret}
% \end{equation}
\begin{align}
    R_1&:= \sum_{t=0}^T   \bigg\{ x^{\intercal}_tP(\theta_t)x_t-\mathds{E}\left[x^{\intercal}_{t+1}P(\theta_{t+1})x_{t+1}|\mathcal{F}_{t}\right]\bigg\}, \notag\\
    R_2&:= \sum_{t=0}^T\mathds{E}\big[x^{\intercal}_t\left(P(\theta_{t+1})-P(\theta_t)\right)x_t|\mathcal{F}_t\big], \label{R2}\\
R_3&:= \sum_{t=0}^T\bigg\{ \left(A\ust x_t+B\ust u_t\right)^{\intercal} P(\theta_t)\left(A\ust x_t+B\ust u_t\right) -\left(A_tx_t+B_tu_t\right)^{\intercal} P(\theta_t)\left(A_tx_t+B_tu_t\right)  \bigg\}, \notag\\ 
    R_4&:= \sum_{t=0}^T\left(J(\theta_t)-J(\theta\ust)\right). \notag
\end{align}
\end{lemma}
%{\color{red} Defined $\cF_t$?}
\begin{proof}
Consider an algorithm that implements $u_t=K(\theta_t)x_t$ at time $t$. Note that 
$x_{t+1}=A\ust x_t+B\ust u_t+w_{t+1}$. Define $\tilde{x}^u_{t+1}:=A_t x_t+B_t u+w_{t+1}$. Then,
the Bellman optimality equation for the Linear Quadratic control problem can be written as follows,
\begin{align*}
J^\star(\theta_t)+x^{\intercal}_t P(\theta_t)x_t&=\min_{u} \bigg(x^{\intercal}_t Q x_t  +u^{\intercal} R u   
 +\mathds{E}\left[{{(\tilde{x}^{u}_{t+1})^{\intercal}}}P(\theta_t)\tilde{x}^u_{t+1}| \mathcal{F}_t\right] \bigg)\\
 &= \bigg(x^{\intercal}_t Q x_t  +u^{\intercal}_t R u_t   
 +\mathds{E}\left[{{(\tilde{x}^{u_t}_{t+1})^{\intercal}}}P(\theta_t)\tilde{x}^{u_t}_{t+1}| \mathcal{F}_t\right] \bigg).
\end{align*}
Upon substituting the value of $\tilde{x}^{u_t}_{t+1}$ in the above, we get
\begin{align}
& J^\star(\theta_t)+x^{\intercal}_t P(\theta_t)x_t \notag\\
& = \left(x^{\intercal}_t Q x_t + u^{\intercal}_t R u_t\right) 
+\mathds{E}\Big[(A_t x_t+B_t u_t+w_{t+1})^{\intercal} P(\theta_t)(A_t x_t+B_t u_t+w_{t+1})| \mathcal{F}_t\Big].\label{eq:rs_4}
\end{align}

%Upon re-arranging the above, we get,
%\begin{align}
%J^\star(\theta_t)+x_t^T P(\theta_t)x_t - \left(x_t^T Q x_t + u_t^T R u_t\right) = 
%\mathds{E}\Big[(A_t x_t+B_t u_t+w_{t+1})^T P(\theta_t)(A_t x_t+B_t u_t+w_{t+1})| \mathcal{F}_t\Big].\label{eq:rs_4}
%\end{align}

Note that $w_{t+1}=x_{t+1}-\left(A\ust x_t+B\ust u_t\right)$ and $w_{t}$ is a martingale difference sequence (Assumption \ref{assumption:noise}). Thus, the l.h.s. of~\eqref{eq:rs_4} can be written as follows,
\begin{align*}
J^\star(\theta_t)+x^{\intercal}_t P(\theta_t)x_t - \left(x^{\intercal}_t Q x_t + u^{\intercal}_t R u_t\right) = &J^\star(\theta_t)-J^\star(\theta\ust)+x^{\intercal}_tP(\theta_t)x_t \\
&- \left(x^{\intercal}_t Q x_t + u^{\intercal}_t R u_t -J^\star(\theta\ust)\right)\\
&=\mathds{E}\left[x^{\intercal}_{t+1}P(\theta_{t+1})x_{t+1}|\mathcal{F}_{t}\right] \\
&- \mathds{E}\big[x^{\intercal}_t\left(P(\theta_{t+1})-P(\theta_t)\right)x_t|\mathcal{F}_t\big]\\&
-  \left(A\ust x_t+B\ust u_t\right)^{\intercal}P(\theta_t)\left(A\ust x_t+B\ust u_t\right)\\ &+\left(A_tx_t+B_tu_t\right)^{\intercal}P(\theta_t)\left(A_tx_t+B_tu_t\right) .
\end{align*}
Therefore by taking a sum from $t=0$ to $t=T$ on both sides, we get $R(T)=R_1+R_2+R_3+R_4$.
\end{proof} 

\begin{lemma}\label{lemma:r1}
On the event $\cE_1 \cap \cE_2$, we have
$
R_1 \leq 2DW^2\sqrt{2T\log \frac{8}{\delta}}+n \sqrt{B_\delta}
$,
where
$B_\delta := b\log \left( \frac{4n\sqrt{b}}{v\delta} \right)$, $b:=v+T(cDX_T)^2(1+c_0^2)$, 
$W=nL\sqrt{2n\log\left(\frac{8nT}{\delta}\right)}$ and $D$ is as in~\eqref{def:D}.
\end{lemma}

\begin{proof}
The proof is the same as that of Lemma 7 in \cite{abbasi_11}, and hence omitted.
\end{proof}
\begin{lemma}\label{lemma:r2}
On the event $\cE_1 \cap \cE_2$,~we have $R_2 \leq2DX_T^2\log_2 T$,
where $X_T$ is defined in \eqref{eq:x}.
\end{lemma}
\begin{proof}
The term in the summation $\sum_{t=0}^T\mathds{E}\big[x_t^{\intercal}\left(P(\theta_{t+1})-P(\theta_t)\right)x_t|\mathcal{F}_t\big]$ (\ref{R2}) corresponding to time $t$ is non-zero only when a change in the policy occurs at $t$. There are $(n + m) \log_2
\left(1 + T X_T^2
\frac{1 + c^{2}_0
}{\lambda}\right)$ episodes till time $T$ (Lemma 8, \cite{abbasi_11}). Therefore, there are $(n + m) \log_2
\left(1 + T X_T^2
\frac{1 + c^{2}_0
}{\lambda}\right)$ non-zero terms and each of them is bounded by $2DX_T^2$.
\end{proof}
\begin{lemma}\label{lemma:r3}
On the event $\cE_1 \cap \cE_2$, we have $
R_3 \leq \frac{8cX_T^2D(1+c_0^2)}{\sqrt{\lambda}}\sqrt{T\beta_T\left(\frac{\delta}{8}\right)\log\frac{det(Z_T)}{det(\lambda I)}}
$,
where $Z_T :=\max_{0 \leq t \leq T}\|z_t\|,$ and $\beta_T\left(\frac{\delta}{8}\right)$ is defined in \eqref{eq:conf}.
\end{lemma}
%\subsection{Proof of Lemma \ref{lemma:r3}}
\begin{proof} The proof is the same as that of Lemma 13 in \cite{abbasi_11}, and hence omitted.
%\rs{The proof is almost the same as that of Lemma 13 in \cite{abbasi_11}, with one minor difference. In our case, as a new episode begins when either $t=2^k$, or $det(Z_t)>2det(Z_0)$, while~\cite{abbasi_11} begins only when $det(Z_t)>2det(Z_0)$. The proof of Lemma 13~\cite{abbasi_11} in turn depends upon Lemma~12\cite{abbasi_11}, where the ratio $\frac{det(Z_t)}{det(Z0)}$ is bounded by $2$ due to the fact that a new episode begins as soon as the determinant doubles. In our case also, it follows from the definition of $\tau_k$, that the ratio $\frac{det(Z_t)}{det(Z0)}$ can be bounded by $2$.} 
\end{proof}
\begin{lemma}\label{thm:r4}
On the event $\cE_1 \cap \cE_2$,
$R_4 \leq \frac{1}{\alpha_0}\left( \beta_{T}\left(\frac{\delta}{4}\right) +\lambda c^2  \right)\sqrt{T}$.

% \begin{equation*}
% R_4 \leq \frac{\sqrt{2}}{\alpha_0(\sqrt{2}-1)}\left( \beta_{T}\left(\frac{\delta}{4}\right) +\lambda c^2  \right)\sqrt{T}.
% \end{equation*}

\end{lemma}
\label{apndix:r4}
\begin{proof}
As defined in Lemma \ref{lemma:decomposition}, we have,
\begin{equation*}
R_4=\sum_{t=0}^T \left(J\ust(\theta_t)-J\ust(\theta\ust)\right).
\end{equation*}
During the $k$-th episode, the algorithm chooses 
$u_t=K(\theta_{t_k})x_t, ~\forall~ t = t_k,t_k +1,\ldots,t_{k+1}$, 
where, $\theta_{t_k}$ is as in~\eqref{eq:rbmle}, and obtained by solving the corresponding optimization problem at the beginning of the episode at time $t_k$. Therefore $R_4$ can be written as :
\begin{equation*}
R_4=\sum_{k=0}^K \Delta_k, ~\text{where}~ \Delta_k:= (t_{k+1}-t_k)\left(J\ust(\theta_{t_k})-J\ust(\theta\ust)\right).
\end{equation*}
$\Delta_k$ is bounded as follows:
\begin{align}
\Delta_k=& (t_{k+1}-t_k)\left( J(\theta_{t_k})-J(\theta\ust) \right)
\leq \frac{(t_{k+1}-t_k)}{\alpha(t_k)}\left(V_{t_k}(\theta\ust)-V_{t_k}(\theta_{t_k})\right),\label{ineq:rs_4}
\end{align}
where the inequality holds since $\theta_{t_k}$ is a minimizer of $V_{t_k}(\theta)+ \alpha(t_k)J^{\star}(\theta)$~\eqref{eq:rbmle}. 
Moreover, 
\begin{align*}
V_{t_k}(\theta\ust)-V_{t_k}(\theta_{t_k}) 
&= V_{t_k}(\theta\ust)+V_{t_k}(\hat{\theta}_{t_k})-V_{t_k}(\hat{\theta}_{t_k})-V_{t_k}(\theta_{t_k})\\
& \leq V_{t_k}(\theta\ust)-V_{t_k}(\hat{\theta}_{t_k}),
\end{align*}
%Since $V_t(\theta) = E_t(\theta)  - \lambda||\theta||^2$, we have,
%\begin{align*}
%V_{t_k}(\theta\ust)-V_{t_k}(\theta_{t_k}) &= (E_{t_k}(\theta\ust)-\lambda||\theta\ust||^2)
%- (E_{t_k}(\theta_{t_k})-\lambda||\theta_{t_k}||^2)\\
%& \leq E_{t_k}(\theta\ust)+E_{t_k}(\hat{\theta}_{t_k})-E_{t_k}(\hat{\theta}_{t_k})-E_{t_k}(\theta_{t_k}) + \lambda||\theta_{t_k}||^2\\
%&\leq E_{t_k}(\theta\ust)-E_{t_k}(\hat{\theta}_{t_k})+\lambda c^2,
%\end{align*}
where the inequality follows since $\hat{\theta}_{t_k}$ is a minimizer of $V_{t_k}(\cdot)$.

Since $\theta\ust \in C_{t_k}(\delta)$, it follows from the definition of the confidence ball that $ V_{t_k}(\theta\ust)-V_{t_k}(\hat{\theta}_{t_k}) \leq \beta_{t_k}\left(\frac{\delta}{4}\right)$. Since $ \beta_{t}(\delta\slash 4) \le \beta_{T}(\delta\slash 4)$,~we have
$V_{t_k}(\theta\ust)-V_{t_k}(\hat{\theta}_{t_k}) \leq\beta_{T}\left(\frac{\delta}{4}\right)$.~Therefore,
\begin{equation}
V_{t_k}(\theta\ust)-V_{t_k}(\theta_{t_k}) \leq \beta_{T}\left(\frac{\delta}{4}\right).    
\end{equation}
Setting $\alpha(t)=a_0\sqrt{T} $, we get 
\begin{align*}
\sum_{k=1}^{K} \Delta_k &\leq \beta_{T}\left(\frac{\delta}{4}\right)\sum_{k=1}^{K} \frac{t_{k+1}-t_k}{\alpha_0\sqrt{T}}= \frac{1}{\alpha_0}  \beta_{T}\left(\frac{\delta}{4}\right)\sqrt{T}.
\end{align*}
\end{proof}
\subsection{Proof of Theorem \ref{thm:r}}
\label{apndix:r}
\begin{proof}
To analyze regret on the event $\cE_1 \cap \cE_2$, we substitute individual bounds on $R_1,~R_2,~R_3$ and $R_4$ in order to obtain
\begin{align*}
R(T) \leq &2DW^2\sqrt{2T\log \frac{8}{\delta}}+n \sqrt{B_\delta}+ 2DX_T^2\log_2 T+ \frac{8X_T^2SD(1+C^2)}{\sqrt{\lambda}}\sqrt{T\beta_T\left(\frac{\delta}{8}\right)\log\frac{det(Z_T)}{det(\lambda I)}} \\
&+  \frac{1}{\alpha_0} \beta_{T}\left(\frac{\delta}{4}\right)\sqrt{T}. 
\end{align*}
\end{proof}
\section{Definition of $d_t$}
The quantity $d_t$ in the definition of $\cE_2(t)$ in (\ref{E2}) is defined as follows,
\begin{align*}
    d_t :&=\frac{1}{1-\rho}\left(\frac{\eta}{\rho}\right)^{n+m}\bigg[2L\sqrt{n\log \frac{4nt(t+1)}{\delta}}+G\mathcal{Z}_T^{\frac{n+m}{n+m+1}}\beta_t\left(\frac{\delta}{4}\right)^{(2(n+d+1))^{-1}}\bigg],\\
    \eta :&= \max \left\{1,\sup_{\theta \in \mathcal{S}}||A\ust+B\ust K(\theta)||\right\},\\
    \mathcal{Z}_T :&=\max_{0 \leq t \leq T}\|z_t\|,\\
    G :&=2 \left(\frac{2c(n+m)^{n+m+0.5}}{\sqrt{U}}\right)^{(n+m+1)^{-1}}, \\%\prk{depends \ on m}
    U :&=\frac{U_0}{H},\\
    U_0 :&=\frac{1}{16^{n+m-2}\max\{1,c^{2(n+m-2)}\}},\\
    H &\text{ is constant such that }
    H >\max\left\{16,\frac{4c^2M^2}{(n+m)U_0}\right\},\\
    M :&= \sup_{Y \geq 0} \frac{nL\sqrt{(n+m)\log \left(\frac{1+\frac{TY}{\lambda}}{\delta}\right)}+\lambda^{1/2}c}{Y}.
\end{align*}

\section{Simulation Experiments}\label{apndix:sims}
In this section, we provide the details on the simulation experiments, along with some additional results. The code and instructions for replicating the presented results are provided in the supplementary material.
%\subsection{Algorithms}
%\subsection{Experimental Setup}
\begin{enumerate}[leftmargin=*]
\item We begin by describing the linear systems used for our experiments in Section \ref{sec:simulations}.
\begin{enumerate}
\item \textbf{Unmanned Aerial Vehicle (UAV): } This system represents a linearized dynamics of an unmanned aerial vehicle (UAV) in a two-dimensional plane, which has been recently studied in the context of reinforcement learning in \cite{lale22,Zhao21}. The  first and third states represent the positions, while the second and fourth states represent the velocities in each dimension. The inputs are accelerations in each dimension. 
\begin{equation*}
    A\ust=\begin{bmatrix}
1 &0.5& 0& 0\\
0&1& 0& 0\\
0& 0& 1& 0.5\\
0& 0& 0&1
\end{bmatrix},
B\ust=\begin{bmatrix} 
0.125& 0\\
0.5& 0\\
0 &0.125\\
0 &0.5
\end{bmatrix},
Q=diag(1,0.1,2,0.2), R=I_2.
\end{equation*}
\item \textbf{Unstable Laplacian Dynamics} This represents a Laplacian system where the adjacent nodes are weakly connected. The lack of stablility (i.e., $\lambda_{max}(A \ust)\geq 1$) makes it a challenging example for system identification and hence it has been studied recently in \cite{dean_17,abbasi_19,dean_20,tu_18,lale22}.  The system matrices are as follows:
\begin{equation*}
    A\ust=\begin{bmatrix}
1.01 & 0.01 & 0\\
0.01 & 1.01 & 0.01\\
0 & 0.01 & 1.01
\end{bmatrix},
 B\ust=I_3,Q=I_3,R=I_3.
\end{equation*}\\
    \item \textbf{Large transient dynamics:}  
We also consider the following unstable system which additionally exhibits large transients.
\begin{equation*}
    A\ust=\begin{bmatrix}
1 & 0 & 0\\
1.1 & 1 & 0\\
0 & 1.1 & 1
\end{bmatrix},
 B\ust=I_3,Q=I_3,R=I_3.
\end{equation*}
\item  \textbf{Longitudinal Flight Control of Boeing 747:} This represents the linearized dynamics of Boeing 747 at 40,000 ft altitute and speed of 774 ft/sec, which was first introduced in \cite{ishihara}. The empirical performance of OFULQ, TS and StabL for this system was recently studied in \cite{lale22}. The four states represent velocity of aircraft along the body axis, velocity perpendicular to the body axis, angle of the body axis with horizontal and the angular velocity. The inputs are elevator angle and thrust of the aircraft. The system matrices are as follows:
    \begin{equation*}
    A\ust=\begin{bmatrix}
0.99 &0.03 &-0.02& -0.32\\
0.01 &0.47 &4.7& 0\\
0.02 &-0.06& 0.4& 0\\
0.01 &-0.04& 0.72& 0.99\\
\end{bmatrix},
 B\ust=\begin{bmatrix} 0.01 &0.99\\
-3.44& 1.66\\
-0.83& 0.44\\
-0.47& 0.25 \\\end{bmatrix}
,Q=I_4,R=I_4 .
\end{equation*}
\end{enumerate}
\item In our experiments, we compared the empirical performance of Augmented RBMLE-UCB and RBMLE with following algorithms: (1) OFULQ \cite{abbasi_11}, (2) Thompson Sampling  \cite{abeille_17}, (3) StabL \cite{lale22}, (4) Randomized Certainty Equivalence (RCE) \cite{faradonbeh_rce}, and (5) Input Perturbations \cite{faradonbeh2018input}. 
The pseudo-code for all of the implemented algorithms is given in Algorithm \ref{algo2}, where the choice of $\theta_{t_k}$ and $u_t$ made by each algorithm are described in Table \ref{table:algo}. The optimization problems for ARBMLE, RBMLE, OFULQ and StabL described in Table \ref{table:algo} are non-convex problems. We used projected gradient descent to solve the optimization problems.~Expression for gradient of the RBMLE objective with respect to $\te$ can be obtained explicitly as in~\cite{abbasi2013online_thesis}.
\begin{algorithm}[h]
   \caption{Reinforcement Learning for LQ systems.}
   \label{algo2}
\begin{algorithmic}
\STATE {\bfseries Initialize:} $t=0$,~$Z_0=\lambda I_{n+m}$
\FOR{$k=0,1,\cdots$}
\IF{$det(Z_t)>2det(Z_{t_{k-1}})$} 
\STATE Calculate $\theta_t$ as defined by the RL algorithm (See Table \ref{table:algo}).
\ELSE 
\STATE $\theta_t=\theta_{t-1}$
     \ENDIF
    \STATE $u_t=f(K(\theta_t),x_t) $ (See Table \ref{table:algo}).
    \STATE $Z_{t+1}=Z_t+z_sz_s^T$
   \ENDFOR
\end{algorithmic}
\end{algorithm}
\begin{center}
\begin{table}[h]
    \centering
\begin{tabular}{|c |c | c |} 
\hline
Algorithm & $\theta_{t_k}$ & $u_t, ~\forall ~t ~\in \left\{t_k,t_k +1,\ldots,t_{k+1}-1\right\}$\\
 \hline
ARBMLE & $ \arg \min_{\theta \in \mathcal{S} \cap \mathcal{C}_{{t_k}}(\delta)} \left\{V_{t_k}(\theta)+ \alpha({t_k})J^{\star}(\theta)\right\},$ & $K(\theta_{t_k})x_t$  \\ \hline
RBMLE & $ \arg \min_{\theta \in \mathcal{S} } \left\{V_{t_k}(\theta)+ \alpha(t)J^{\star}(\theta)\right\},$ & $K(\theta_{t_k})x_t$ \\ \hline
OFULQ & $\arg \min_{\theta \in \mathcal{S} \cap \mathcal{C}_{{t_k}}(\delta)} J^{\star}(\theta)$ & $K(\theta_{t_k})x_t$ \\ \hline
TS & $\hat{\theta}_{t_k}+\beta_{t_k}(\delta)Z_{t_k}^{-1/2}\mathcal{N}(0,1)$ & $K(\theta_{t_k})x_t$ \\ \hline
IP & $\hat{\theta}_{t_k}$ & $K(\theta_{t_k})x_t + \eta^{IP}_t$  \\ \hline
RCE & $\hat{\theta}_{t_k} +\eta^{RCE}_t$ & $K(\theta_{t_k})x_t$  \\ \hline
StabL & $\arg \min_{\theta \in \mathcal{S} \cap \mathcal{C}_{{t_k}}(\delta)} J^{\star}(\theta)$ & $\begin{cases}
K(\theta_{t_k})x_t + \mathcal{N}(0,\sigma_w^2) \text{ if }  T<T_w \\
K(\theta_{t_k})x_t  \text{ otherwise }
\end{cases}$  \\ \hline
\end{tabular}
    \caption{Choices of $\theta_{t_k}$ and $u_t$ for various algorithms.}
    \label{table:algo}
\end{table}
\end{center}
%\item
%First begin by describing qualitatively the types of systems you have simulated.
%Do also prominently mention that many of the system are open-loop unstable.
%You can describe what type of system each is, i.e., its characteristics. For example, one describes Boeing, or ..., etc. Also, where the system came from (for. e.g., which paper), or if it was randomly generated. Why you chose that system to test.
%Also, provide the actual equations of the systems, including all the values of the system parameters
%(e.g., noise mean and variance, distribution, initial condition of the system, matrices, cost matrices, etc.).
%\item
%Next, discuss the nature of all the policies that were simulated.
%You can qualitatively describe each, as well as provide the control laws.
%Table \ref{table:algo} would be good here.
\item Initially, the controls are chosen as follows in order to obtain a initial estimates of the system: $$u_t=K_{\text{init}}x_t + \eta_t\text{ for } 0\leq t\leq T_{\text{init}} \text{ and } \eta_t \text{ is }  \mathcal{N}(0,1), $$ i.i.d., $K_{\text{init}}$ is a stabilizing controller and $x_{0}=0$. The noise is pre-generated, ensuring that initialization is uniform across algorithms.
The definition of confidence interval in ARBMLE, TS and StabL depends on the choice of confidence parameter $\delta$ and a constant $c$ such that $trace({\theta^\star}^T\theta^\star) <c^2$. StabL algorithm uses an excitation $\mathcal{N}(0,\sigma_w^2)$ for $T<T_w$.  The values of various hyper-parameters used in our experiments are described in Table \ref{table:hyper}.
\begin{center}
\begin{table}[h]
    \centering
\begin{tabular}{|c|c|c|c|c|c|c|c|} 
\hline
Parameter & $T_{\text{init}}$ & $T$ & $\delta$ & $\lambda$ & $\alpha_0$ & $\sigma_w$ & $T_w$  \\ \hline
Value & 50 & 500 & $10^{-4}$ & $10^{-4}$ & $10^{-2}$ & 2 & 35  \\\hline
\end{tabular}
    \caption{Values of various parameters}
    \label{table:hyper}
\end{table}
\end{center}
%\item
%Next discuss the role and character of all the hyperparameters in the textual description. Qualitatively.
%\item
%Next describe how you initialize the algorithms. Mention the role of $T_{\init}, K_{init}$, etc. Just qualitatively.
%\item
%Then create a Table of Hyperparameters vs Systems, containing VALUES of everything: $T_{init}, K_{init}, \delta$, etc., for each system for each algorithm. 
%Also discuss any issues involved with each parameter choice.
\item
We provide simulation results for following additional examples. Figure \ref{fig:a} includes the comparison between ARBMLE, OFULQ, TS and StabL. Figure \ref{fig:b} provides a comparison between ARBMLE,~RBMLE,~RCE and IP. The average regret values for these systems at $T=500$ are shown in Table \ref{table:3}.
\begin{figure}[h]
     \centering
     \begin{subfigure}{0.45\textwidth}
         \centering
     \includegraphics[width=\textwidth]{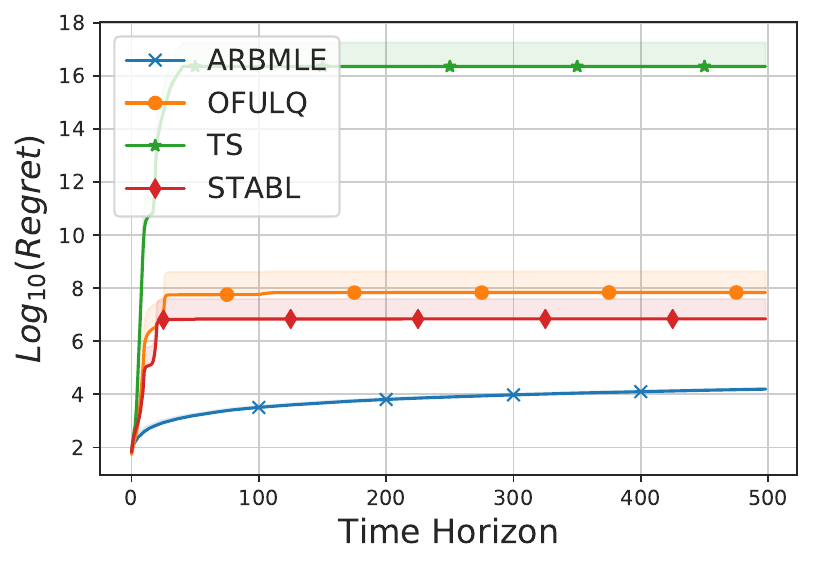}
         \caption{Stabilizable, not controllable}
         \label{fig:a_1}
     \end{subfigure}
     \begin{subfigure}{0.45\textwidth}
         \centering
     \includegraphics[width=\textwidth]{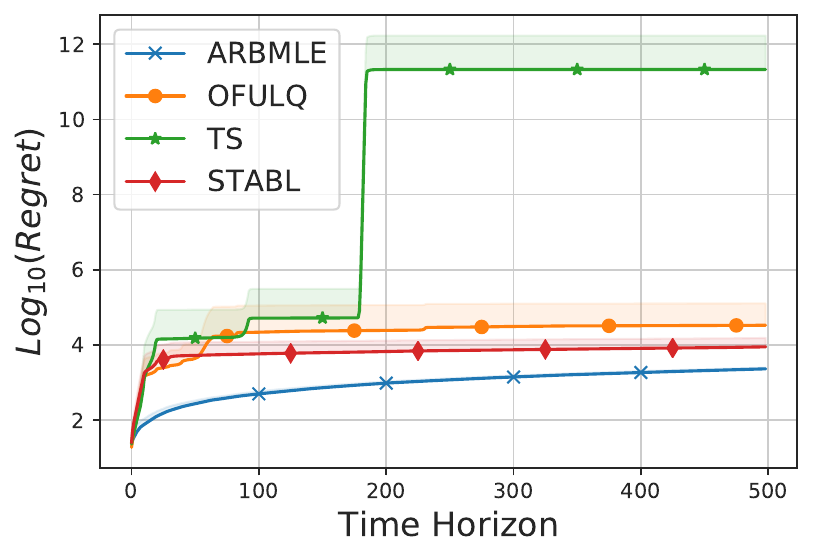}
         \caption{Chained Integrator Dynamics }
              \label{fig:a_2}
    \end{subfigure}
    \caption{Logarithm of the Averaged Regret over 50 runs of  ARBMLE, OFULQ, TS. and StabL.}
        \label{fig:a}
\end{figure}
\begin{figure}[h]
     \centering
     \begin{subfigure}{0.45\textwidth}
         \centering
     \includegraphics[width=\textwidth]{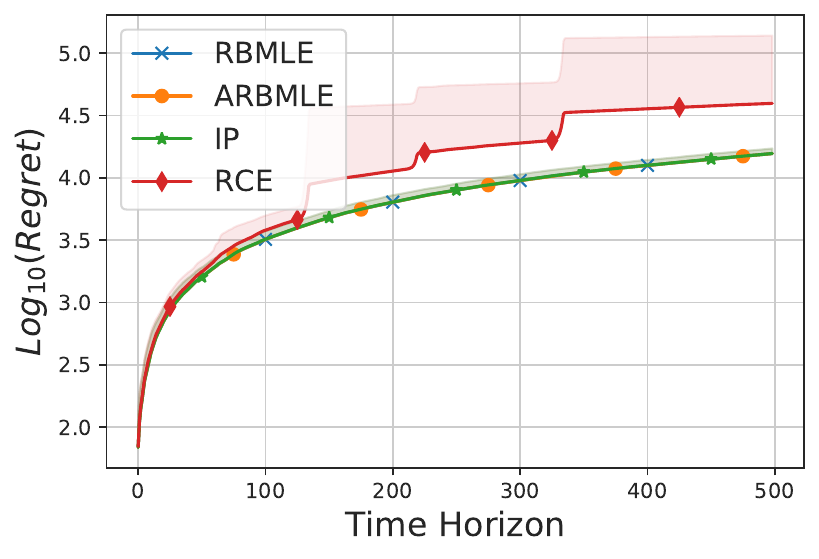}
         \caption{Stabelizable, not controllable}
         \label{fig:b_1}
     \end{subfigure}
     \begin{subfigure}{0.45\textwidth}
         \centering
     \includegraphics[width=\textwidth]{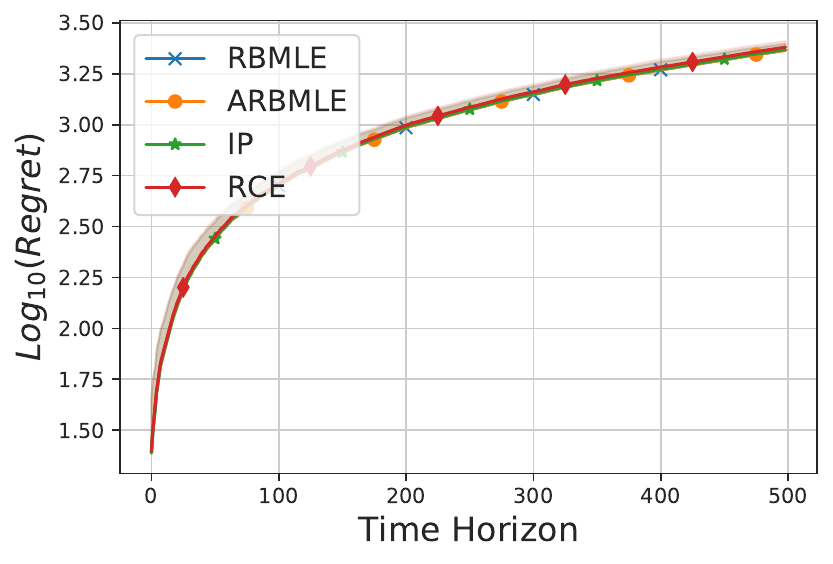}
         \caption{Chained Integrator Dynamics }
              \label{fig:b_2}
    \end{subfigure}
    \caption{Logarithm of the Averaged Regret over 50 runs of ARBMLE, RBMLE, RCE, TS }
        \label{fig:b}
\end{figure}
\item 
\begin{enumerate}
\item Stabilizable but Not Controllable System:
We consider a system studied in \cite{lale22} which is stabilizable but not controllable. Lack of controllability is challenging for system identification. ARBMLE/RBMLE outperforms OFULQ, TS, StabL and RCE by a significant margin. 
\begin{equation*}
 A\ust=\begin{bmatrix}
-2& 0 &1.1\\
1.5& 0.9& 1.3\\
0 &0& 0.5\\
\end{bmatrix},
 B\ust=\begin{bmatrix}
 1 &0\\
0 &1\\
0& 0\\
\end{bmatrix}
,Q=R=I_3 .
\end{equation*}

\item Chained Integrator Dynamics:
We consider a simple chained integrator system with 2-dimensional states and 2-dimensional input. 
\begin{equation*}
    A\ust=\begin{bmatrix}
1 & 0.1 \\
0 & 1 \\
\end{bmatrix},
 B\ust=I_2,Q=I_2,R=I_2.
\end{equation*}
\begin{comment}

    \item Non-Sparse Matrix: We consider the following non-sparse example with $n=m=4$ studied in \cite{farjadnasab2022model}.
\begin{equation*}
    A\ust=\begin{bmatrix}
0.3 & 0.4 & 0.2 & 0.2\\
0.2 & 0.3 & 0.2 & 0.2\\
0.2 & 0.2 & 0.4 & 0.4\\
0.4 & 0.2 & 0.2 & 0.4\\
\end{bmatrix},
 B\ust=I_4,Q=R=I_4 .
\end{equation*}
\end{comment}
\begin{center}
\begin{table}[h]
    \centering
\begin{tabular}{|c |c | c | c| c|c| c|c|} 
\hline
Ex. & RBMLE & ARBMLE & OFULQ & TS & IP & RCE & STABL\\
 \hline
(a) & 15665 & 15663 &  $6.9 \times 10^{7} $&  $2.2 \times 10^{16} $  &  15628 & 39593 & $6.9 \times 10^{6}$ \\ \hline
(b)  &2322 & 2322 & 33449 & $2.1 \times 10^{11}$  & 2337 & 2402 &8927\\ \hline
%(c)& 16144& 16135 & $2.1 \times 10^{12}$ & $2.1 \times 10^{11}$ &16164  & $2.1 \times 10^{11}$ & $1,2 \times 10^{9}$  \\ \hline
\end{tabular}
    \caption{Average Regret Performance at $T=500$.}
    \label{table:3}
\end{table}
\end{center}
\end{enumerate}

\item
Additonal Remarks:
\begin{itemize}[leftmargin=1pt]
    \item We use a stabilizing controller $K_{\text{init}}$ for initialization of our simulation experiments similar to \cite{dean_17}. Note that our theoretical regret analysis does not assume knowledge of a stabilizing controller, unlike some recent works on adaptive control of LQG systems including \cite{dean_17,faradonbeh_rce,simchowitz20}.
    \item As demonstrated in the simulation results, OFULQ and Thompson Sampling have a~very large initial regret indicating poor initial estimates of system parameters (also highlighted in \cite{lale22}). ARBMLE/RBMLE, IP and RCE show much better initial regret performance compared to OFULQ and TS.     
    \item Implementation of ARBMLE, OFULQ, StabL and TS involve definition $\beta_t(\delta) $ which denotes boundary of confidence interval. Our simulations for ARBMLE, OFULQ, StabL and TS are based on $\beta_t(\delta) $ as defined in \eqref{eq:conf}.
Instead, recent works \cite{dean_17,lale22} use  $\beta_t(\delta) := trace\left((\theta\ust-\hat{\theta}_t)^{\intercal}Z_t(\theta\ust-\hat{\theta}_t)\right)$.  However, one may note that $\theta\ust$ in $\beta_t(\delta)$ is not known to the learning agent, and so such a definition of $\beta_t(\delta)$ is not a viable for implementation. The effect of the choice of $\beta_t(\delta)$ on the regret performance is shown in the Figure \ref{fig:conf}.
\begin{figure}[t]
     \centering
     \begin{subfigure}{0.32\textwidth}
         \centering
     \includegraphics[width=\textwidth]{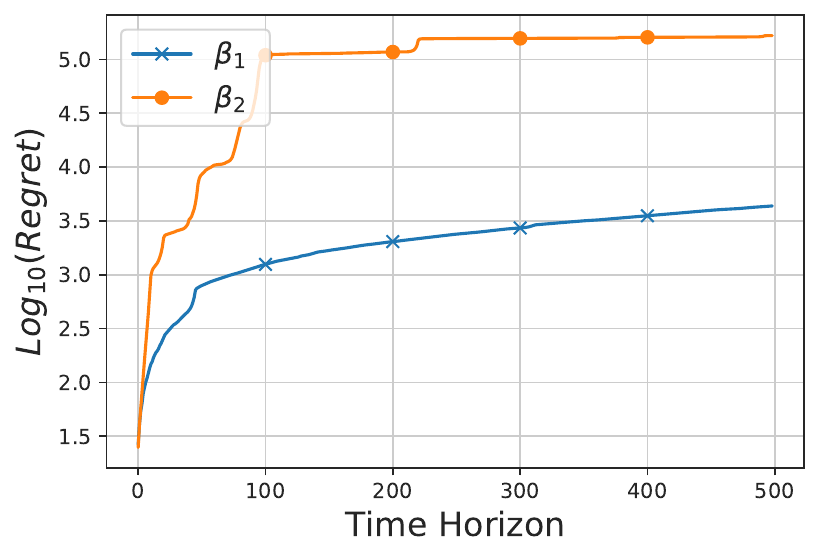}
         \caption{OFULQ}
         \label{fig:conf_1}
     \end{subfigure}
     \begin{subfigure}{0.32\textwidth}
         \centering
     \includegraphics[width=\textwidth]{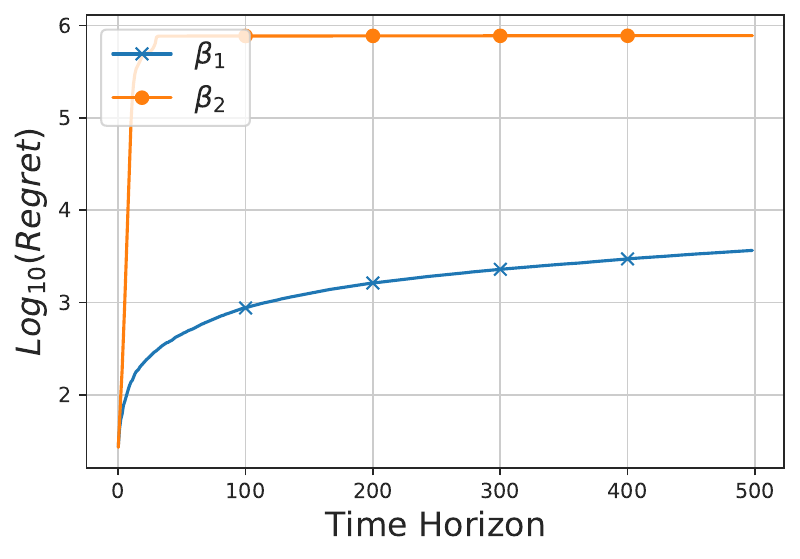}
         \caption{TS }
              \label{fig:conf_2}
    \end{subfigure}
     \begin{subfigure}{0.32\textwidth}
         \centering
     \includegraphics[width=\textwidth]{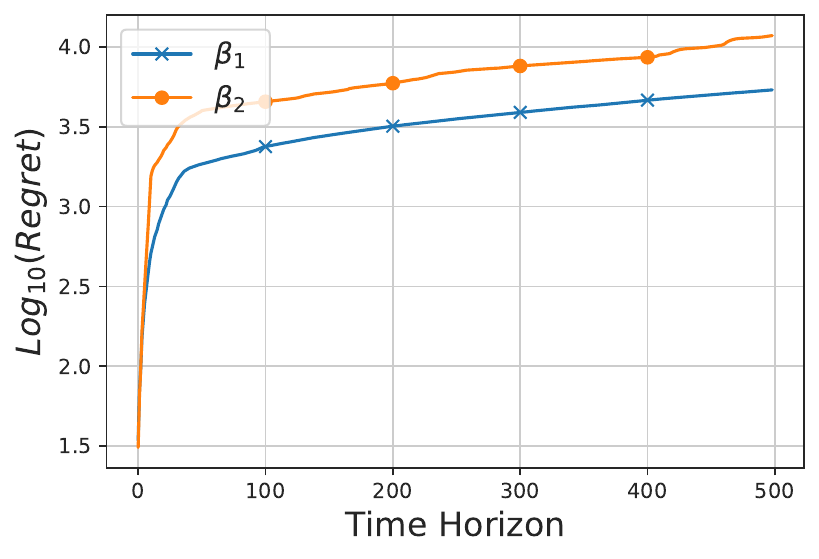}
         \caption{Stabl}
              \label{fig:conf_3}
    \end{subfigure}
    \caption{Effect of choice of confidence interval definition on performance. $\beta_1$: confidence interval as defined in \ref{eq:conf}. $\beta_2$: confidence interval as defined in \cite{dean_17}   }
     \label{fig:conf}
\end{figure}
\begin{comment}
\item  The pre-constant in the bias-term $\alpha_0$ balances the trade-off between exploration and exploitation, specifically whether to choose a $\theta_t$ with lower average cost or a $\theta_t$ closer to the least-squared estimate. In Figure \ref{fig:bias}, we demonstrate the effect of choice of $\alpha_0$ on the regret performance.
\begin{figure}[h]
     \centering
     \includegraphics[width=0.5\textwidth]{finalplots/bias (2).pdf}
    \caption{Performance of ARBMLE for different values of $\alpha_0$.}
 \label{fig:bias}
\end{figure}
\end{comment}

\item 
The estimates of OFULQ lies on boundary, while the estimates of ARBMLE\slash RBMLE, IP and RCE are closer to the least squared estimate. Note that RBMLE can be seen as Lagrangian version of OFULQ, indicating that $ \alpha(t)$ may be much smaller than the implicit Lagrange multiplier for OFULQ.
\end{itemize}
%\item The code for our simulation experiments is provided in the supplementary material. The seed values for random number generation are set appropriately for replication of the results. The instruction for the code are provided  in supplementary material.
%\item
%Indicate how the code can be accessed, how the random numbers can be accessed.
%Provide all the information needed for the results to be replicated. And say that you have provided all the information for the results to be replicated.
\end{enumerate}

\end{document}